\documentclass[11pt]{amsart}
\usepackage{eucal}
\usepackage{tikz}
\usepackage{xypic}
\usepackage{amsfonts,amssymb}
\usepackage{epsfig}
\usepackage{enumerate}
\usepackage{enumitem}
\usetikzlibrary{decorations.pathreplacing}

\newtheorem{theorem}{Theorem}[section]
\newtheorem{maintheorem}[theorem]{Main Theorem}
\newtheorem{lemma}[theorem]{Lemma}
\newtheorem{prop}[theorem]{Proposition}

\newtheorem{corollary}[theorem]{Corollary}
\newtheorem{defn}[theorem]{Definition}

\newcommand{\Q}{\mathbb Q}
\newcommand{\C}{\mathbb C}
\newcommand{\R}{\mathbb R}
\newcommand{\K}{\mathbb K}
\renewcommand{\P}{\mathbb P}

\newcommand{\D}{\mathbb D}

\newcommand{\Z}{\mathbb Z}
\newcommand{\N}{\mathbb N}

\newlist{HF}{enumerate}{1}
\setlist[HF]{label=(HF\arabic*)}

\newlist{HF6}{enumerate}{1}
\setlist[HF6]{label=(HF6)}

\newlist{HF7}{enumerate}{1}
\setlist[HF7]{label=(HF7)}

\newlist{PB}{enumerate}{1}
\setlist[PB]{label=(PB\arabic*)}

\newlist{gi}{enumerate}{1}
\setlist[gi]{label=(gi\arabic*)}

\newlist{NS}{enumerate}{1}
\setlist[NS]{label=(NS\arabic*)}

\newlist{gfun}{enumerate}{1}
\setlist[gfun]{label=(g\arabic*)}

\title[Growth rate of wrapped Floer cohomology]{Affine varieties, Singularities and the Growth rate of wrapped Floer cohomology}
\author{Mark McLean}

\begin{document}

\begin{abstract}
In this paper, we give partial answers to the following questions:
Which contact manifolds are contactomorphic to
links of isolated complex singularities?
Which symplectic manifolds are symplectomorphic to smooth affine varieties?
The invariant that we will use to distinguish such manifolds
is called the growth rate of wrapped Floer cohomology.

Using this invariant we show that if $Q$ is a simply connected manifold whose unit cotangent bundle is contactomorphic to the link of an isolated singularity or whose cotangent bundle is symplectomorphic
to a smooth affine variety
then M must be rationally elliptic and so it must have certain bounds on its Betti numbers.
\end{abstract}

\maketitle

\bibliographystyle{halpha}


\tableofcontents

\section{Introduction}

To any affine variety $B \subset \C^N$ with an isolated
singularity at $0$, we can form its {\it link}
which is a $C^\infty$ manifold given by
the intersection of $B$ with a small sphere.
The link of any isolated singularity has a natural contact structure
(\cite{Varchenko:isolatedsingularities}).
If $A \subset \C^N$ is a smooth affine variety then
it has a natural symplectic structure given by restricting
the standard symplectic structure on $\C^N$ to $A$.

In this paper, we give partial answers to the following questions:
Which contact manifolds are contactomorphic to
links of isolated complex singularities?
Which symplectic manifolds are symplectomorphic to smooth affine varieties?
The invariant that we will use to distinguish such manifolds
is called the growth rate of wrapped Floer cohomology.
\begin{theorem} \label{theorem:showtheorem}
Let $Q$ be a compact oriented Riemannian manifold.
Suppose that the unit cotangent bundle of $Q$
is contactomorphic to the link of an isolated
complex singularity or that $T^*Q$ is symplectomorphic to
a smooth affine variety.
\begin{enumerate}
\item If $Q$ is a simply connected $4$ or $5$ manifold
then it must be homeomorphic to one of
$S^4,\C \P^2, S^2 \times S^2,\C \P^2 \# \overline{\C \P^2},\C \P^2 \# \C \P^2$
or diffeomorphic to one of
$S^5,S^3 \times S^2,SU(3)/SO(3)$ or the non-trivial
$S^3$ bundle over $S^2$.
\item
If $Q$ is simply connected of any dimension $n$ then its
$j$th Betti number is bounded above by
$\frac{1}{2} \left (
\begin{array}{c}
   n \\
   j
\end{array}
\right ).$
\item The fundamental group of $Q$
satisfies certain growth constraints
(explained below).
\end{enumerate}
\end{theorem}

We will now state the main theorem of this paper
and explain why it implies the above theorem.
The contact manifolds that we will be looking at will be boundaries of certain symplectic manifolds
called {\it Liouville domains} (Definition \ref{defn:liouvilledomain}).
The symplectic manifolds that we will be examining are
obtained from Liouville domains by attaching a
cylindrical end and extending the symplectic form.
For a Liouville domain $M$, such a symplectic manifold
is called the {\it completion}
$\widehat{M}$ of $M$.
Let $L_1,L_2$ be two oriented
exact Lagrangians in $\widehat{M}$ which are cylindrical outside $M$ (Definition \ref{defn:admissiblelagrangian})
admitting a spin structure.
From now on we will call such objects {\it admissible Lagrangians}.
Then for any field $\K$, we can define an invariant $\Gamma(L_1,L_2;\K) \in  \{-\infty\} \cup [0,\infty]$
called the {\it growth rate of wrapped Floer cohomology}
(Definition \ref{defn:growthrate}).


\begin{maintheorem} \label{theorem:main}
Suppose that $\partial M$ is contactomorphic to
the link of an isolated complex singularity
or $\widehat{M}$ is symplectomorphic to a
smooth affine variety.
Then $\Gamma(L_1,L_2;\K) \leq n$  for any
transversally intersecting admissible Lagrangians $L_1,L_2$ in $\widehat{M}$ and any field $\K$ where $n$ is the complex dimension of
our variety.
\end{maintheorem}

We give a sketch of the proof of this Theorem
in Section \ref{section:sketch}.
Let $Q$ be a compact oriented Riemannian manifold
and let $q \in Q$ be a basepoint.
The {\it based loopspace of length $\leq \lambda$} for $\lambda \in [0,\infty]$ written as
$\Omega^{\leq \lambda} Q$ is defined to be the space of continuous maps $\R / \Z \to Q$ sending $0$ to $q$
of length $\leq \lambda$ equipped with the $C^0$ topology.
We define $\Omega(Q) := \Omega^{\leq \infty}(Q)$.
For any field $\K$, define $g(Q,\lambda;\K)$ to be the rank of the image of the natural map
$H_*(\Omega^{\leq \lambda};\K) \rightarrow H_*(\Omega(Q);\K)$ and
define $\Gamma(Q;\K) := \limsup\limits_{\lambda\rightarrow \infty}\frac{\log{g(Q,\lambda;\K)}}{\log{\lambda}}$.
Such a number does not depend on the choice of basepoint or Riemannian metric on $Q$.

\begin{corollary} \label{corollary:main}
Suppose $Q$ is a Riemannian manifold satisfying at least one of the following conditions:
\begin{enumerate}[label=(\alph*)]
\item The unit cotangent bundle $S^* Q$ is contactomorphic to the link of some isolated singularity or
\item \label{item:affinesympcondition}
$T^* Q$ is symplectomorphic to an affine variety.
\end{enumerate}
Then $\Gamma(Q;\K) \leq \text{dim}_\R Q$ for every field $\K$.
\end{corollary}

This will follow from the results
in \cite{AbbondandoloSchwartz:conormalboundary}.
The proof is contained in subsection
\ref{subsection:proofofcorollarymain}.
The condition that  $\Gamma(Q;\K) \leq \text{dim}_\R Q$ for every field $\K$ is a very restrictive condition.
Here are some constraints:
\begin{enumerate}
\item If $Q$ is an oriented surface then it must have genus $\leq 1$.
\item \label{item:gromovresult}
If $Q$ has finite fundamental group then
by \cite{Gromov:homotopical},
we get that
\[\limsup_i \left(\frac{\log(\sum_{j \leq i} \text{dim}(H_j(\Omega(Q);\K)))}{\log(i)}\right)  \leq \text{dim}_\R Q\]
for all $\K$.
In particular the sum of the first $i$ Betti numbers of the
based loopspace of $Q$
is bounded above by a polynomial in $i$.
Hence, if $Q$ is simply connected then
$Q$ is rationally elliptic
(see \cite{FelixHalperinThomoasRadiusofConvergence}).
Rationally elliptic means that
$\text{dim}(\pi_*(Q) \otimes \Q) < \infty$.

\item If $Q$ is a simply connected $4$ manifold
then it must be homeomorphic to one of
\[S^4,\C \P^2, S^2 \times S^2,\C \P^2 \# \overline{\C \P^2}, \C \P^2 \# \C
\P^2\] (\cite[Lemma 3.2]{PaternainPetean:minimalentropy} combined with (\ref{item:gromovresult})).
\item If it is a simply connected $5$ manifold then it must be diffeomorphic to
$S^5,S^3 \times S^2,SU(3)/SO(3)$ or the non-trivial
$S^3$ bundle over $S^2$
(\cite[Corollary 3.6]{PaternainPetean:minimalentropy} combined with (\ref{item:gromovresult})).
\item
If $Q$ is simply connected of any dimension $n$ then its Betti numbers $b_i$ satisfy:
\[b_i \leq \frac{1}{2} \left (
\begin{array}{c}
   n \\
   i
\end{array}
\right ).\]
See \cite[Corollary to Theorem 1]{Pavlov:bettiestimates}
combined with (\ref{item:gromovresult}).
\item The fundamental group of $Q$ has growth rate at most $\text{dim}_\R Q$.
Here the growth rate of a finitely generated group $G$ is defined
as follows: choose a finite set of generators and let
$l(\lambda)$ be the number of elements of $G$ expressed in words of these generators of length $\leq \lambda$,
then
the {\it growth rate of} $G$ is defined as $\limsup\limits_{\lambda\rightarrow \infty} \frac{\log{l(\lambda})}{\log{\lambda}}$.
\end{enumerate}
Theorem \ref{theorem:showtheorem}
follows immediately.

Corollary \ref{corollary:main} part \ref{item:affinesympcondition} has a relationship with the following conjecture.
We will suppose for this conjecture that $Q$ is simply connected.
We say that $Q$ has a {\it good complexification} if there exists
a smooth affine variety $U$ defined over $\R$ diffeomorphic to $Q$ such that the natural map $U(\R) \rightarrow U(\C)$
is a homotopy equivalence.
There is a question in \cite{Totaro:complexifications}
which asks if $Q$ has a Riemannian metric of non-negative sectional curvature
when $Q$ has a good complexification.
There is also a conjecture attributed to Bott which
says that any simply connected manifold with non-negative sectional curvature is
rationally elliptic (see
\cite[Question 12, page 519]{FelixHalperinThomas:rationalhomotopytheory}).
Felix and Thomas in \cite{FelixThomasRadius} proved that if a simply connected manifold is rationally elliptic
then the growth rate of the rational Betti numbers of its based loopspace grows sub exponentially.
Finally the result in \cite{Gromov:homotopical} tells us that if the Betti numbers grow sub exponentially
for $Q$ then $\Gamma(Q;\Q) < \infty$. Putting all of this together we get the following question:
if $U(\R) \hookrightarrow U(\C)$ is a homotopy equivalence then is $\Gamma(Q,\Q) <\infty$?

This project was partially supported by the National Science Foundation (DMS-1508207).
The author would like to thank Paul Seidel,
Kenji Fukaya and Burt Totaro for helpful remarks.

\section{Basic Definitions and Properties} \label{section:floercohomologydefinition}

\subsection{Liouville Domains and Floer Cohomology} \label{section:liouvilefloerdefn}

\begin{defn} \label{defn:liouvilledomain}
Let $M$ be a compact manifold with boundary and a $1$-form $\theta_M$ satisfying:
\begin{enumerate}
\item $\omega_M := d\theta_M$ is a symplectic form.
\item The $\omega_M$-dual $X_{\theta_M}$ of $\theta_M$ points outwards along $\partial M$.
\end{enumerate}
Then we say that $(M,\theta_M)$ is a {\bf Liouville domain}
and $\theta_M$ is called the associated {\bf Liouville form}.
Because $X_{\theta_M}$ is transverse to $\partial M$,
we get that $\alpha_M := \theta_M|_{\partial M}$
is a contact form.
The {\bf contact boundary} of $M$
is the contact manifold $(\partial M,\text{ker}(\alpha_M))$.
By flowing $\partial M$ backwards along $X_{\theta_M}$ we get that
a small collar neighborhood of $\partial M$ is equal to $(1-\epsilon,1] \times \partial M$
with $\theta_M = r_M \alpha_M$ where $r_M$ parameterizes the interval.
The {\bf completion} $\widehat{M}$ of $M$ is given by
extending this collar neighborhood by
attaching $[1,\infty) \times \partial M$
to $M$ and extending $\theta_M$ by $r_M \alpha_M$.
\end{defn}

\begin{defn} \label{defn:admissiblelagrangian}
	A (possibly non-compact) properly embedded submanifold of $L \subset \widehat{M}$ is said to be an {\bf exact Lagrangian which is cylindrical outside $M$} if
	\begin{itemize}
		\item it is of half the dimension of $\widehat{M}$,
		\item $\theta_M|_{L} = df_L$ for some smooth $f_L : L \to \R$ where $f_L = 0$ outside $M$ and
		\item the vector field $\frac{\partial}{\partial r_M}$ is tangent to $L$ in the cylindrical end
		$[1,\infty) \times \partial M$.
	\end{itemize}
	We say that $f_L$ is a {\bf function associated} to $L$.
	
	An {\bf admissible Lagrangian}
	is an oriented exact Lagrangian which is cylindrical outside $M$ with a choice of spin structure.
\end{defn}

Here are some important examples of Liouville domains
and admissible Lagrangians:

{\it Example 1}: Let $Q$ be a compact Riemannian manifold.
Then we have a canonical $1$ form $\theta_Q$
on $T^*Q$ defined by $\theta_Q(V) := P_*(V)(\beta)$
for all $1$-forms $\beta$ on $Q$ and vectors $V \in T_\beta T^*Q$
where $P : T^*Q \to Q$ is the natural projection.
The bundle of covectors $D^*Q$ of length $\leq 1$
is our Liouville domain with Liouville form $\theta_Q$
and $\widehat{D^*Q}$ is symplectomorphic to $T^*Q$.
This symplectomorphism is
the identity on $D^*Q$ and it sends
a point $(r,s) \in [1,\infty) \times \partial D^* Q$
to $rs \in T^*Q$.
Because the cotangent fiber $T_q^*Q \subset T^*Q$ is
contractible and because
$\theta_Q|_{T_q^* Q} = 0$,
we get that it
is an admissible Lagrangian inside $\widehat{M}$.

{\it Example 2}: Let $A \subset \C^N$ be a smooth
affine variety. This has a natural symplectic structure
$\omega_A$ given by restricting the standard symplectic form
on $\C^N$.
Let $\theta_A = \sum_{j=1}^N \frac{1}{2}r_j^2 d\vartheta_j|_A$
where $(r_j,\vartheta_j)$
are standard polar coordinates
on the $j$th $\C$ factor.
Then the intersection of $A$ with a very large
closed ball in $\C^N$ is a Liouville domain
with Liouville form $\theta_A$ whose
completion is symplectomorphic to $(A,\omega_A)$
(see \cite[Lemma 2.1]{McLean:affinegrowth}).
An example of an admissible Lagrangian inside
an affine variety would be $\R^N \subset \C^N$.

\begin{defn}
A time dependent Hamiltonian $H : S^1 \times \widehat{M} \to \R$ is said to be {\bf admissible} if
it is equal $\lambda r_M + C$ near infinity for some $\lambda, C \in \R$.
The constant $\lambda$ is called its {\bf slope}.

Sometimes we will view $H$ as a time dependent family
of Hamiltonians
$(H_t)_{t \in \R / \Z}$.
\end{defn}
Admissible Hamiltonians on
$T^*Q$ of slope $\lambda$ consist of smooth functions
equal to
$\lambda |\beta|_Q$ for all $\beta \in T^*Q$ outside a compact set
where $|\cdot|_Q$ is the norm induced by our Riemannian metric on $Q$.
The time $1$ flowlines of such Hamiltonians inside a level set near infinity correspond to length $\lambda$ geodesics after projecting to $Q$.
An example of an admissible Hamiltonian
on $\C^n$ is $\lambda |z|^2$
where $|\cdot|$ is the standard Euclidean metric.

We will  give a very brief definition
of Floer cohomology
$HF^*_{[a,b]}(L_1,L_2,H_t)$ for any admissible
Lagrangians $L_1,L_2$ and any admissible Hamiltonian
$(H_t)_{t \in \R / \Z}$
inside the completion $\widehat{M}$
of a Liouville domain $(M,\theta_M)$
(see \cite{Floer:lagrangianmorsecomplex},
\cite{Floer:gradientflow},
\cite{ohfloer1993},
\cite{ohfloer1995addendum}).
The only difference is that the Lagrangians
that we use are non-compact
but this is not a problem
as we can use the maximum principle
\cite[Lemma 7.2]{SeidelAbouzaid:viterbo}.
Floer cohomology for such non-compact Lagrangians
was defined in \cite{SeidelAbouzaid:viterbo},
but we are only interested in Floer cohomology
as a group without additional $A_\infty$ products which means that our definition is more straightforward.

Having said that in this paper
we won't deal with the definition directly
and instead we use properties
\ref{item:upperbound}-\ref{item:translationinvariance}
below (with the exception of the proof of
Corollary \ref{corollary:main} where we use
the additional properties
\ref{item:actionlimit} and
 \ref{item:additionalproperty}
 combined with the fact that we need to
 define $HF^*$ for a slightly more general
 class of Hamiltonians).
From now on we fix a coefficient field $\K$.

Define $\phi_t^{H_t}$ to be the time $t$ flow
of the associated Hamiltonian vector field $X_{H_t}$
defined by $i_{X_{H_t}}\omega_M = -dH_t$.
We will first define this Floer cohomology group
when $\phi^{H_t}_1(L_1)$ and $L_2$ intersect transversally,
and then later on we will define it in the case
where they may not intersect transversally.
We suppose $f_{L_1},f_{L_2}$ are the functions associated to $L_1$ and $L_2$ respectively.
To any intersection point $p \in \phi_1^{H_t}(L_1) \cap L_2$ we have a value called its {\bf action}
\[{\mathcal A}(p) := f_{L_2}(p) - f_{L_1}((\phi_1^{H_t})^{-1}(p))  + \]
\[\int_0^1
 H_t(\phi^{H_t}_{t}(\phi^{H_t}_{-1}(p))) - \theta_M(\frac{d}{dt}\phi^{H_t}_{t}((\phi^{H_t}_1)^{-1}(p))) dt.\]
The chain complex, written as $C^*_{[a,b]}(L_1,L_2,H_t)$, is the free $\K$ vector space generated by
intersection points $\phi^{H_t}_1(L_1) \cap L_2$ whose action is in $[a,b]$.
We can ensure that each intersection point $p$ has an index $|p| \in \Z/2\Z$
which makes $C^*_{[a,b]}(L_1,L_2,H_t)$ into a $\Z/2\Z$-graded $\K$ vector space.
In this paper we are not really concerned with gradings
and so we will not define them here.
The differential depends on a choice of almost complex structure but
the resulting cohomology group does not depend on this.
We say that $J$ is a {\bf cylindrical almost complex structure}
on $M$ if it is compatible with the symplectic form $\omega_M$
(i.e. $\omega_M(\cdot,J(\cdot))$ is a Riemannian metric)
and if $\theta_M \circ J = dr_M$ outside a large compact set.
Choose a smooth $S^1$ family of cylindrical almost
complex structures $(J_t)_{t \in \R/\Z}$. 
We define ${\mathcal M}(p,q,H_t,J_t)$ to be the set of smooth maps
\[ u : \R \times [0,1] \to M \]
satisfying:
\[ \partial_s u + J_t \partial_t u = J_t X_{H_t},\]
\[u(s,0) \in L_1, \quad u(s,1) \in L_2,\]
\[ u(s,t) \to p(t) ~ \text{as} ~ s \to -\infty \quad
\quad \text{and} \quad u(s,t) \to q(t) ~  \text{as} ~ s \to +\infty.\]
There is a free $\R$
action on ${\mathcal M}(p,q,H_t,J_t)$
given by translation in the $s$ coordinate.
For generic $J_t$
and for $|p|=|q|+1$,
$\overline{{\mathcal M}}(p,q,H_t,J_t):= {\mathcal M}(p,q,H_t,J_t)/\R$
is a disjoint union of manifolds $\sqcup_j \overline{{\mathcal M}}^j(p,q,H_t,J_t)$
where $\overline{{\mathcal M}}^j(p,q,H_t,J_t)$ has dimension $j$.
Also the spin structure ensures that
$\overline{{\mathcal M}}^0(p,q,H_t,J_t)$ is an oriented $0$ dimensional manifold.
We define
$\#\overline{{\mathcal M}}^0(p,q,H_t,J_t)$
 for $|p|=|q|+1$
as the number of elements in
$\overline{{\mathcal M}}^0(p,q,H_t,J_t)$
counted with sign.
The differential is linear and on intersection points $q \in \phi^{H_t}_1(L_1) \cap L_2$ it is defined as:
\[ \partial(q) := \sum_{|p| = |q|+1}\#
\overline{{\mathcal M}}^0(p,q,H_t,J_t) p.\]
The differential increases the action,
and so this chain complex has a natural filtration
by action $p \to {\mathcal A}(p)$.
We define $HF^*_{[a,b]}(L_1,L_2,H_t)$ to be the homology of the above chain complex.

If $H^1_t,H^2_2$ are two admissible Hamiltonians
with $H^1_t \leq H^2_t$ for all $t \in \R / \Z$,
we have a natural map:
$HF^*_{[a,b]}(L_1,L_2,H^1_t) \to HF^*_{[a,b]}(L_1,L_2,H^2_t)$.
This is called a {\bf continuation map}
and it is defined in a similar way
by counting certain moduli spaces.

Now suppose that $\phi_1^{H_t}(L_1)$ does not intersect $L_2$ transversally and that $a,b$
is not in the image of ${\mathcal A}$. Then we define
$HF^*_{[a,b]}(L_1,L_2,H_t)$ to be the direct limit
\[\varinjlim_{H'_t} HF^*_{[a,b]}(L_1,L_2,H'_t)\]
where
$H'_t < H_t$ are admissible Hamiltonians
so that $\phi^1_{H'_t}(L_1)$ and $L_2$ intersect transversally and the directed system
is taken with respect to the ordering $\leq$.
Sometimes one cannot find such admissible
Hamiltonians $H'_t$ which $C^\infty$
converge to $H_t$.
In this case one needs to have
more general Hamiltonians (see \cite[Section 8a]{SeidelAbouzaid:viterbo}).
Note that we can define $HF^*_{(a,b]}(L_1,L_2,H_t)$, $HF^*_{[a,b)}(L_1,L_2,H_t)$ and
$HF^*_{(a,b)}(L_1,L_2,H_t)$ in a similar way.
If $a = -\infty$ and $b = +\infty$,
we will write $HF^*(L_1,L_2,H_t) := HF^*_{[a,b]}(L_1,L_2,H_t)$.
Also we will define
$HF^*_-(L_1,L_2,H_t) := HF^*_{(-\infty,0]}(L_1,L_2,H_t)$.

\subsection{Properties of Floer Cohomology} \label{section:propertiesofgrowthrate}

\begin{HF}

\item \label{item:upperbound}
The rank of $HF^*_{[a,b]}(L_1,L_2,H_t)$ is bounded above by the number
of intersection points $\phi^{H_t}_1(L_1) \cap L_2$ whose action is in $[a,b]$
assuming that all such intersection points
are transversal.
\item \label{item:actionmorphismproperty}
If $a_1 \geq a_2$ and $b_1 \geq b_2$ then there is a natural morphism:
\[ HF^*_{[a_1,b_1]}(L_1,L_2,H_t) \to HF^*_{[a_2,b_2]}(L_1,L_2,H_t).\]
We call such a morphism an {\bf action morphism}.
Composing two action morphisms gives another action morphism.
Similar properties hold for other intervals
of the form $[a,b)$, $(a,b]$ and $(a,b)$.
Such a morphism is an isomorphism if there are no intersection points
$p \in \phi^{H_t}_1(L_1) \cap L_2$
of action ${\mathcal A}(p)$ in the interval $[a_2,b_1] \setminus [a_1,b_2]$.
For $-\infty \leq a \leq b \leq c \leq \infty$ we have the following long exact sequence:
\[ \rightarrow HF^*_{(b,c]} (L_1,L_2,H_t) \rightarrow HF^*_{[a,c]} (L_1,L_2,H_t) \rightarrow HF^*_{[a,b]} (L_1,L_2,H_t) \rightarrow \]
where the morphisms $HF^*_{(b,c]} (L_1,L_2,H_t) \rightarrow HF^*_{[a,c]} (L_1,L_2,H_t)$
and $HF^*_{[a,c]} (L_1,L_2,H_t) \rightarrow HF^*_{[a,b]} (L_1,L_2,H_t)$ are action morphisms.
\item \label{item:continuationmorphismproperty}
If $H_{1,t} \leq H_{2,t}$ then there is a natural morphism
\[ HF^*_{[a,b]}(H_{1,t}) \to HF^*_{[a,b]}(H_{2,t}).\]
This is a {\bf continuation morphism}
and the composition of two such morphisms is also a continuation morphism.
\item \label{item:commutativityproperty}
We have that continuation morphisms commute with action morphisms as follows:
Let $a_1,a_2,b_1,b_2 \in \R$ so that $a_1 \geq a_2$, $b_1 \geq b_2$ and let $H_{1,t} \leq H_{2,t}$ be Hamiltonians.
Then we have the following commutative diagram:
\[
\xy
(0,0)*{}="A"; (60,0)*{}="B";
(0,-20)*{}="C"; (60,-20)*{}="D";
"A" *{HF^*_{[a_1,b_1]}(L_1,L_2,H_{1,t})};
"B" *{HF^*_{[a_2,b_2]}(L_1,L_2,H_{1,t})};
"C" *{HF^*_{[a_1,b_1]}(L_1,L_2,H_{2,t})};
"D" *{HF^*_{[a_2,b_2]}(L_1,L_2,H_{2,t})};
{\ar@{->} "A"+(18,0)*{};"B"-(18,0)*{}};
{\ar@{->} "C"+(18,0)*{};"D"-(18,0)*{}};
{\ar@{->} "A"+(0,-4)*{};"C"+(0,4)*{}};
{\ar@{->} "B"+(0,-4)*{};"D"+(0,4)*{}};
%
\endxy
\]
where the horizontal arrows are action morphisms
and the vertical arrows are continuation morphisms.
\item \label{item:translationinvariance}
We have an isomorphism $HF^*(L_1,L_2,H_t) \to HF^*(L_1,L_2,H_t+c)$.
If $c>0$ then this isomorphism is induced
by the natural continuation map.
If $c<0$ then it is induced by the inverse
of a continuation map.
%
\end{HF}

\subsection{Growth Rate of Wrapped Floer Cohomology} \label{section:growthratedefn}

All Floer cohomology groups are taken with coefficients in a field $\K$ unless stated otherwise.
Let $(M,\theta_M)$ be a Liouville domain.
Let $L_1,L_2$ be admissible Lagrangians inside $\widehat{M}$.
We let $H \geq 0$ be an admissible Hamiltonian with positive slope.
By property \ref{item:continuationmorphismproperty}
we have natural continuation maps
$HF^*(L_1,L_2,\lambda_1 H) \to HF^*(L_1,L_2, \lambda_2 H)$
for $\lambda_1 \leq \lambda_2$.
\begin{defn} \label{defn:growthrate}
	Define the {\bf wrapped Floer cohomology group}
	of $L_1$ and $L_2$ to be 
	\[HW^*(L_1,L_2,H,\K) := \varinjlim_\lambda HF^*(L_1,L_2,\lambda H).\]
	The {\bf growth rate} $\Gamma(L_1,L_2,H;\K) \in \{-\infty\} \cup [0,\infty]$ is defined as:
	\[ \limsup_\lambda \frac{\log{a_\lambda}}{\log{\lambda}} \]
	where $a_\lambda$ is the rank of the image of the natural map:
	\[ HF^*(L_1,L_2,\lambda H) \to HW^*(L_1,L_2,H,\K).\]
	If the rank $a_\lambda$ is zero then we define $\log{a_\lambda} := -\infty$.
	
	Lemma \ref{lemma:independentofh} below tells
	us that $\Gamma(L_1,L_2,H;\K)$ does not depend on $H$.
	Hence we will write
	$\Gamma(L_1,L_2;\K)$ for the growth rate.
	We will sometimes write $\Gamma(L_1,L_2)$ or $\Gamma(L_1,L_2,H)$ when its is clear which coefficient field we are using.
\end{defn}

The following two Lemmas are used to show us that $\Gamma(L_1,L_2,H;\K)$
does not depend on the choice
of Hamiltonian $H$. Lemma \ref{lemma:hamiltonianswithsameslope}
will also be used later on.

\begin{lemma} \label{lemma:hamiltonianswithsameslope}
	Let $H_1 \leq H_2$ be two admissible Hamiltonians
	with the same slope and let $L_1,L_2$
	be admissible Lagrangians.
	Then the continuation map
	\[HF^*(L_1,L_2,H_1) \to HF^*(L_1,L_2,H_2)\]
	is an isomorphism.
\end{lemma}
\proof
Because $H_1,H_2$ have the same slope,
there is a constant $C>0$ so that
$H_1 + C > H_2$.
By \ref{item:continuationmorphismproperty}
this means that we get the following sequence
of continuation maps:
\begin{multline} \label{eqn:fourarrows}
	HF^*(L_1,L_2,H_2 - C) \to 
	HF^*(L_1,L_2,H_1) \xrightarrow{e} \\
	HF^*(L_1,L_2,H_2) \to
	HF^*(L_1,L_2,H_1+C).
\end{multline}
By \ref{item:translationinvariance}
and the fact that the composition of any two continuation
maps is a continuation map by \ref{item:continuationmorphismproperty},
the composition of any two continuation maps from
(\ref{eqn:fourarrows}) is an isomorphism.
Hence $e$ is an isomorphism.
\qed

\begin{lemma} \label{lemma:independentofh}
	Let $H_1,H_2 \geq 0$ be admissible Hamiltonians
	and let $L_1,L_2$ be admissible Lagrangians.
	Then:
	\[\Gamma(L_1,L_2,H_1;\K) = \Gamma(L_1,L_2,H_2;\K).\]
\end{lemma}
\proof
Let $\kappa_1 > 0$ be the slope of $H_1$
and let $\kappa_2 >0$ be the slope of $H_2$.
Choose an admissible Hamiltonian
$H_3 \geq 0$ with slope equal to $\kappa_1$
so that
$H_3 \geq H_1$ and $H_3 \geq \frac{\kappa_1}{\kappa_2}H_2$.
By Lemma \ref{lemma:hamiltonianswithsameslope},
the continuation maps
\[HF^*(L_1,L_2,H_1) \to HF^*(L_1,L_2,H_3)\]
and
\[HF^*(L_1,L_2,\frac{\kappa_1}{\kappa_2}H_2) \to HF^*(L_1,L_2,H_3)\]
are isomorphisms.
Also \ref{item:continuationmorphismproperty}
tells us that these two continuation maps
commute with all other continuation maps.
Hence
\[\Gamma(L_1,L_2,H_1) = \Gamma(L_1,L_2,H_3)
= \Gamma(L_1,L_2,\frac{\kappa_1}{\kappa_2}H_2).\]
Also $\Gamma(L_1,L_2,\frac{\kappa_1}{\kappa_2}H_2)$
is equal to
$\Gamma(L_1,L_2,H_2)$. Hence
$\Gamma(L_1,L_2,H_1;\K) = \Gamma(L_1,L_2,H_2;\K).$

\qed

\subsection{Proof of Corollary \ref{corollary:main}}
\label{subsection:proofofcorollarymain}

Part of the proof
of Corollary \ref{corollary:main}
needs the additional properties
\ref{item:actionlimit} and
\ref{item:additionalproperty} below for $HF^*$
combined with the fact that we need to
define $HF^*$ for a slightly more general
class of Hamiltonians.
This is only contained in the proof of the preliminary
technical Lemma \ref{lemma:loopspacegrowthhfgrowth}.
We need some other definitions and Lemmas before we
prove Corollary \ref{corollary:main}.
The purpose of these Lemmas is to
translate the conventions in \cite{AbbondandoloSchwartz:conormalboundary}
to the conventions in this paper.
From now on we fix our field $\K$.
Every Floer group will be defined over this field.

\begin{defn}
	Let $P : [0,\infty) \to [0,\infty)$ be any function.
	We define $\Gamma(P) := \limsup_\lambda 
	\frac{\log{P(\lambda)}}{\log{\lambda}}$.
	Here $\log(0)$ is defined to be $-\infty$.
	
	If we have a sequence $(p_i)_{i \in \N}$
	then we define $\Gamma((p_i)_{i \in \N}) := \Gamma(P_{\lfloor \rfloor})$
	where
	$P_{\lfloor \rfloor} : [0,\infty) \to [0,\infty)$ satisfies
$P_{\lfloor \rfloor}(x) = p_{\lfloor x \rfloor}$ where $\lfloor x \rfloor$ is the largest integer
		less than or equal to $x$.
\end{defn}

The above definition will be used throughout this paper.
We will also use the fact that if
$P$ is a non-decreasing function
then
$\Gamma(P) = \Gamma((p_i)_{i \in \N})$
where $p_i := P(i)$ for all $i \in \N$.

\begin{defn} \label{defn:sequencegrowthrate}
 Let $(M,\theta_M)$ be a Liouville domain,
 let $(H_i : \widehat{M} \to \R)_{i \in \N}$
 be a sequence of admissible Hamiltonians
 and let $\nu : \N \to (0,\infty]$ be a non-decreasing function so that
 the slope of $H_i$ is $i$ times the slope of $H_1$ and
 $H_i \leq H_{i+1}$ for all $i \in \N$.
 Then we have natural maps:
 \[\iota_{i_1,i_2} : HF^*_{[-\nu(i_1),\infty)}(L_1,L_2,H_{i_1}) \to 
 HF^*_{[-\nu(i_2),\infty)}(L_1,L_2,H_{i_2})\]
 for every $i_1 \leq i_2$ given by the composition of
 the action morphism 
 \[HF^*_{[-\nu(i_1),\infty)}(L_1,L_2,H_{i_1}) \to 
 HF^*_{[-\nu(i_2),\infty)}(L_1,L_2,H_{i_1})\]
 and the continuation morphism
  \[HF^*_{[-\nu(i_2),\infty)}(L_1,L_2,H_{i_1}) \to 
  HF^*_{[-\nu(i_2),\infty)}(L_1,L_2,H_{i_2}).\]
 
 By \ref{item:commutativityproperty},
 we have that for all $i_1 \leq i_2 \leq i_3$,
 $\iota_{i_2,i_3} \circ \iota_{i_1,i_2} = \iota_{i_1,i_3}$ and therefore these maps
 form a directed system.
 Hence we define:
 \[HW^*(L_1,L_2,(H_i)_{i \in \N},\nu) :=
 \varinjlim_i HF^*_{[-\nu(i),\infty)}(L_1,L_2,H_i)\]
 and
 $\Gamma(L_1,L_2,(H_i)_{i \in \N},\nu)) :=
 \Gamma((b_i)_{i \in \N})$
 where
 \[b_i := \text{rank image}(HF^*_{[-\nu(i),\infty)}(L_1,L_2,H_i) \to
 HW^*(L_1,L_2,(H_i)_{i \in \N},\nu)).\]
\end{defn}

\begin{lemma} \label{lemma:growthrateequivalencedefn}
	Let $(H_i)_{i \in \N}$ and $\nu$ be as in Definition
	\ref{defn:sequencegrowthrate} as above.
	Suppose that for all $i \in \N$,
	all intersection points
	$\phi^{H_i}_1(L_1) \cap L_2$ have action in
	$[-\nu(i),\infty)$.
	Then:
	$\Gamma(L_1,L_2) = \Gamma(L_1,L_2,(H_i)_{i \in \N},\nu)$.
\end{lemma}
\proof
Because every intersection point
$\phi^{H_i}_1(L_1) \cap L_2$
has action in $[-\nu(i),\infty)$,
we get by \ref{item:actionmorphismproperty} that the natural action morphism:
\[HF^*_{[-\nu(i),\infty)}(L_1,L_2,H_i)
\to HF^*(L_1,L_2,H_i)\]
is also an isomorphism.
Using this fact combined with the fact that
action morphisms and continuation morphisms
commute by \ref{item:continuationmorphismproperty} and \ref{item:commutativityproperty}, we get:
\[\Gamma(L_1,L_2,(H_i)_{i \in \N},\nu)
= \Gamma(L_1,L_2,(H_i)_{i \in \N},\infty).\]
So from now on, we may as well assume that
$\nu(i) = \infty$ for all $i \in \N$.

Let $H \geq 0$ be an admissible Hamiltonian
with slope equal to the slope of $H_1$.
Because
\[A(\lambda) := \text{rank image} (HF^*(L_1,L_2,\lambda H) \to HW^*(L_1,L_2,H))\]
is non-decreasing in $\lambda$,
we get:
$ \Gamma(A) = \Gamma((a_i)_{i \in \N})$
where
$a_i := A(i)$ for all $i$.

Because the slope of $H_i$ is $i$ times the slope
of $H_1$, we can choose an increasing sequence of constants
$(C_i)_{i \in \N}$
so that
$i H + C_i \geq H_i$ for all $i \in \N$.
By \ref{item:translationinvariance},
we get:
$\Gamma((a_i)_{i \in \N})
= \Gamma(L_1,L_2,(i H + C_i)_{i \in \N},\infty)$.

Also by Lemma \ref{lemma:hamiltonianswithsameslope},
the
natural continuation morphism
$HF^*(L_1,L_2,H_i) \to
HF^*(L_1,L_2,i H + C_i)$
is an isomorphism.
Hence because continuation morphisms commute by \ref{item:continuationmorphismproperty},
$\Gamma(L_1,L_2,(H_i)_{i \in \N},\infty)
= \Gamma(L_1,L_2,(H + C_i)_{i \in \N},\infty)$.
Putting everything together we get:
$\Gamma(L_1,L_2)
= \Gamma(L_1,L_2,(H_i)_{i \in \N},\nu)$.
\qed

\bigskip

Example 1 in Section \ref{defn:liouvilledomain} tells us
that for any choice of metric on $Q$,
\begin{itemize}
	\item $T^* Q$ is the completion of a Liouville domain,
	\item the boundary of this Liouville domain is contactomorphic to the unit cotangent bundle of $Q$
	\item and the fibers $T_q^* Q$ are admissible Lagrangians for all $q \in Q$.
\end{itemize}

\begin{lemma} \label{lemma:loopspacegrowthhfgrowth}
Let $Q$ be an  oriented Riemannian manifold.
Then we have
$\Gamma(T^*_q Q,T^*_p Q;\K) = \Gamma(Q;\K)\quad \forall q,p \in Q$ for some choice of metric on $Q$.
\end{lemma}
\proof
Fix $p,q \in Q$.
Choose a generic metric on $Q$
so that the geodesics joining $p$ and $q$ are
non-degenerate
(I.e. the differential of the exponential map
along each such geodesic is an isomorphism).
Let $|\cdot|$ be the induced norm on $T^* Q$.

The radial coordinate
 on the cylindrical end
is given by $r_Q : T^* Q \to \R$,
$r_Q(\beta) = |\beta|$.
We define $H := r_Q^2$.
Let $L_1 := T^*_q Q$ and $L_2 := T^*_p Q$
be our admissible Lagrangians.
Because all geodesics joining $p$ and $q$
are non-degenerate, we get that
$\phi_1^H(L_1)$ and $L_2$ intersect transversally.
Even though $H$ is not admissible, we can
still define $HF^*_{[a,b]}(L_1,L_2,H)$ in the usual way
as described in Subsection \ref{section:liouvilefloerdefn}.
In fact these Floer cohomology groups can be defined
for any Hamiltonian which is a non-decreasing function
of $|\beta|$ near infinity.
One can also define action morphisms and continuation maps
between such Hamiltonians as well so that they satisfy the same properties.

Define:
\[ b(\lambda) := \text{rank image}(HF^*_{[-\lambda^2,\infty)}(L_1,L_2,H)
\to HF^*(L_1,L_2,H)).\]
For all $c$, all intersection points $\phi^H_1(L_1) \cap L_2$
contained inside $H^{-1}(c)$
have action $-c^2$.
Hence the main result in \cite{AbbondandoloSchwartz:conormalboundary}
tells us that
$\Gamma(Q;\K) = \Gamma(b(\lambda))$.
We have the following property:
\begin{HF6} 
	\item \label{item:actionlimit}
	The natural map
	$\varinjlim_\lambda HF^*_{[-\lambda^2,\infty)}(L_1,L_2,H)) \to
	HF^*(L_1,L_2,H)$
	is an isomorphism.
\end{HF6}
Hence:
\begin{equation} \label{eqn:directlimitaction}
b(\lambda) = \text{rank image}(HF^*_{[-\lambda^2,\infty)}(L_1,L_2)
\to \varinjlim_l HF^*_{[-l^2,\infty)}(L_1,L_2,H)).
\end{equation}
Because all of the geodesic flowlines joining
$p$ and $q$ are non-degenerate,
we get that the set of action values
${\mathcal S} \subset \R$ of $\phi^H_1(L_1) \cap L_2$
is discrete.
Hence there is a constant $\mu > 0$
so that $i \mu \notin S$ for all $i \in \N$.

Define $b_i := b(i \mu)$ for all $i \in \N$.
Because $b(\lambda)$ is non-decreasing,
we get that
$\Gamma(b(\lambda)) = \Gamma((b_i)_{i \in \N})$
and hence
$\Gamma(Q;\K) = \Gamma((b_i)_{i \in \N})$.

Choose smooth functions
$(f_i : [0,\infty) \to \R)_{i \in \N}$ so that:

\begin{itemize}
	\item \label{item:derivatiesoffl}
	$f_i,f'_i,f''_i \geq 0$,
	\item \label{item:constantderivativeatinfintyoffl}
	$f'_i(x) = i\mu + \epsilon_i$ for $x$ sufficiently large where $\epsilon_i > 0$
	is small enough so that $(i\mu,i\mu + \epsilon_i) \cap {\mathcal S} = \emptyset$.
	\item \label{item:mainlyxsq}
	If $x \leq \frac{1}{2}(i\mu)^2$
	then $f_\lambda(x) = x^2$.
\end{itemize}

We define $H_i : T^*Q \to \R$ by $H_i(\beta) =f_i(|\beta|)$.
Now the additional technical property that we
need our Floer cohomology groups to
satisfy is the following:
\begin{HF7}
	\item \label{item:additionalproperty}
	Suppose that we have two functions $f,g : \R \to \R$
	so that $f,f',f''$, $g,g',g'' \geq 0$, $f \leq g$
	and so that
	for all $x \leq S$, $f(x) = g(x)$ for some $S > 0$.
	We will also assume that $F,G : T^*Q \to \R$
	defined by $F(\beta) := f(|\beta|)$
	and $G(\beta) := g(|\beta|)$ are smooth.
	
	Then the natural continuation map
	\[HF^*_{[A,\infty)}(L_1,L_2,F) \to HF^*_{[A,\infty)}(L_1,L_2,G))\]
	is an isomorphism for all $S - Sf'(S) \leq A \leq 0$.
	Such continuation maps also commute with
	other continuation maps and action morphisms
	as well.
\end{HF7}

This property follows from the
maximum principle \cite[Lemma 7.2]{SeidelAbouzaid:viterbo} and the fact that
for any $s \in \R$, the action
of any intersection point
$\phi^H_1(L_1) \cap L_2$ inside $F^{-1}(s)$ is
$sf'(s) - s$ and similarly for $G$.

So by \ref{item:additionalproperty} we have that
the continuation map
\[HF^*_{[-(i\mu)^2,\infty)}(L_1,L_2,H_i) \to HF^*_{[(-i\mu)^2,\infty)}(L_1,L_2,H)\]
is an isomorphism for all $i \in \N$.
Hence by equation (\ref{eqn:directlimitaction})
combined with the fact that
action morphisms and continuation morphisms commute
we get:
\[ \Gamma(b(\lambda)) = \Gamma((b_i)_{i \in \N}) 
= \Gamma(L_1,L_2,(H_i)_{i \in \N},\mu).\]
And so by Lemma
\ref{lemma:growthrateequivalencedefn},
$\Gamma(b(\lambda)) =
\Gamma(L_1,L_2).$
Hence $\Gamma(T^*_q Q,T^*_p Q;\K) = \Gamma(Q;\K)$.
\qed

\bigskip

\proof of Corollary \ref{corollary:main}.
By Lemma \ref{lemma:loopspacegrowthhfgrowth}
we have $\Gamma(T^*_q Q,T^*_p Q;\K) = \Gamma(Q;\K)$ for all $q,p \in Q$.
By Theorem \ref{theorem:main}
we have for all $p \neq q$ that
\[\Gamma(T^*_q Q,T^*_p Q;\K) \leq \text{dim}_\R Q = \text{ the complex dimension of our variety}.\]
Hence $\Gamma(Q;\K) \leq \text{dim}_\R Q$.
\qed

\section{Sketch of the Proof of Main Theorem \ref{theorem:main}} \label{section:sketch}

For simplicity, we will first sketch the proof
in the case when the boundary of our Liouville
domain
is contactomorphic to
the link of a cone singularity.
Most of the key ideas are contained in
the proof of this special case.
We will explain at the end of this section
what needs to be done more generally.

The statement of Main Theorem \ref{theorem:main}
in our special case is:
{\it Let $M$ be a Liouville domain whose contact boundary
is contactomorphic to the link of a cone singularity.
Then for any transversally intersecting
admissible Lagrangians $L_1,L_2$ inside $\widehat{M}$,
we have $\Gamma(L_1,L_2;\K) \leq \frac{1}{2}\text{dim}(M)$.}

The proof of this Theorem splits up in to
three main steps:

{\it Step 1}: We show that if, in some sense,
we can bound the number of flowlines
of $X_H$ joining $L_1$ and $L_2$
for some appropriate Hamiltonian $H$ then
we get a bound on $\Gamma(L_1,L_2;\K)$.
This is the only place where we use the axioms from
Section \ref{section:propertiesofgrowthrate}.
This Step is performed in Section \ref{section:growthratebounds}.

{\it Step 2}: We show that if there is a
relatively compact open set
$U \subset \widehat{M}$ containing $M$
which looks like an a annulus bundle
near $\partial U$ then there is a bound
on the number of flowlines as in Step 1
(see Definition \ref{defn:nicelystratified}
for a more general definition
involving many annulus bundles).
This step is performed in Section \ref{section:nicelystratifiedsection}.

{\it Step 3}: This is the only step where we use the fact
that $\partial M$ is contactomorphic to the link
of a singularity.
We use this contactomorphism
to construct our region $U$
from Step 2.
Hence by Step 2 and Step 1 we get our result.
This step is performed in Section \ref{section:proofofmaintheorem}.

\bigskip

{\bf Step 1}:

The aim of this step is to show that
if $H$ is a Hamiltonian and
$P : [0,\infty) \to [0,\infty)$ a non-decreasing function
so that
\begin{itemize}
\item
$H^{-1}((-\infty,C])$ is a compact
for some $C$,
\item the number of intersection
points $\cup_{\lambda_1 \leq \lambda} (\phi^{H}_{\lambda_1} (L_1) \cap L_2)$
inside each level set of $H$ near $C$
is bounded above by $P(\lambda)$
\item and some other conditions are satisfied
(See Theorem \ref{theorem:growthratebound}),
\end{itemize}
then
$\Gamma(L_1,L_2;\K) \leq \Gamma(P)$.
Note that such a Hamiltonian may not be admissible.

We start this Step by proving the following
{\bf key claim:} 
{\it If there are two admissible Hamiltonians $0 \leq Q_1 \leq Q_2$ equal to
$f_1(r_M)$, $f_2(r_M)$ respectively (pictured below)
where $r_M$ is the cylindrical coordinate
and where the rank of the image of
$HF^*_-(L_1,L_2;\lambda Q_1) \to HF^*_-(L_1,L_2;\lambda Q_2)$
is bounded above by $P(\lambda)$ for all $\lambda$
then $\Gamma(L_1,L_2;\K) \leq \Gamma(P)$}.
(See Lemma \ref{lemma:positiveactionhamiltonian}
for a precise statement).

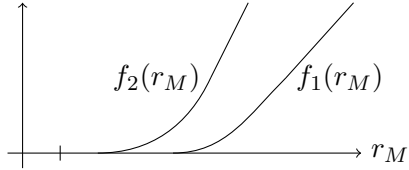
\begin{center}
	\begin{tikzpicture}[domain=0:4]
	\draw[->] (-0.2,0) -- (4.5,0) node[right] {$r_M$};
	\draw[->] (0,-0.2) -- (0,2);
	%
	\draw (0.5,0.1) -- (0.5,-0.1) node[below,scale=1] {};
	%
	\draw (2,0) to[out=0,in=225] (3.5,1) node[right] {$f_1(r_M)$};
	\draw[-] (3.5,1) -- (4.4,2);
	%
	\draw (1,0) to[out=0,in=245] (2.5,1) node[left] {$f_2(r_M)$};
	\draw[-] (2.5,1) -- (3,2);

	\end{tikzpicture}
\end{center}

{\it Sketch of the proof of the key claim:}
One can show that
\[\text{rank image}(HF^*_-(L_1,L_2;\lambda Q_1) \to HF^*_-(L_1,L_2;\lambda Q_2))\]
and \[b_\lambda := \text{rank image}(HF^*(L_1,L_2;\lambda Q_1) \to HF^*(L_1,L_2;\lambda Q_2))\]
differ by a constant by \ref{item:upperbound}
and \ref{item:actionmorphismproperty} due 
	to the fact that
most points in $\phi^{Q_i}_\lambda(L_1) \cap L_2$,$i=1,2$ have negative action.
Because an appropriate multiple $\mu Q_1$ of $Q_1$ has the same slope as
$Q_2$, we get that $b_\lambda$ is equal to
\[c_\lambda := \text{rank image}(HF^*(L _1,L_2;\lambda Q_1) \to HF^*(L_1,L_2;\lambda \mu Q_1))\]
by Lemma \ref{lemma:hamiltonianswithsameslope}.
Also $c_\lambda$  (and hence $b_\lambda$) is bounded below by:
\[\text{rank image}(HF^*(L_1,L_2;\lambda Q_1) \to HW^*(L_1,L_2;Q_1;\K))\]
by \ref{item:continuationmorphismproperty}.
Our key claim now follows from this lower bound on $b_\lambda$.
\qed

\bigskip

Using the above key claim we can complete
Step 1 as follows:
We start with a Hamiltonian $H$ as described above.
We can find appropriate functions
$(g_i : \R \to \R)_{i \in \N}$ (pictured below)
so that the number of intersection
points $\phi^{H_i}_1(L_1) \cap L_2$
is bounded above by $P(\lambda_i)$
where $H_i := g_i \circ H$ and $\lambda_i$
is approximately linear in $i$.

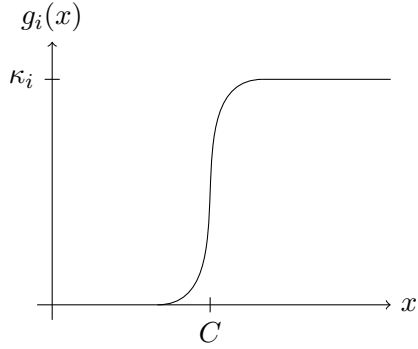
\begin{center}
\begin{tikzpicture}[domain=0:4]
\draw[->] (-0.2,0) -- (4.5,0) node[right] {$x$};
\draw[->] (0,-0.2) -- (0,3.5) node[above] {$g_i(x)$};
\draw (1.4,0) to[out=0,in=180] (2.8,3);
\draw (2.8,3) to (4.5,3);
%
   \draw (0.1,3) -- (-0.1,3) node[left] {$\kappa_i$};
   \draw (2.1,0.1) -- (2.1,-0.1) node[below,scale=1] {$C$};
\end{tikzpicture}
\end{center}
Here $\kappa_i$ is linear in $i$
and the support of the derivative $g'_i(x)$
is a very small neighborhood of $C$.

We let $K$ be an admissible Hamiltonian with small derivatives whose slope is small and positive and which
is $0$ on a very large set.
We can choose such a $K$
so that no additional intersection points $\phi^{H_i + iK}_1(L_1) \cap L_2$
of non-positive action are created.
Hence by \ref{item:upperbound},
one gets that
$HF^*_-(L_1,L_2,H_i + iK) = HF^*_-(L_1,L_2,H_i)$ is bounded above by $P(\lambda_i)$.
One can find Hamiltonians
$0 \leq Q_1 \leq Q_2$ from
the key claim above
so that
$i Q_1 \leq H_i + iK \leq i Q_2$.
This means that we have continuation morphisms
\[
HF_-^*(L_1,L_2;i Q_1) \to
HF_-^*(L_1,L_2;H_i + iK) \to HF_-^*(L_1,L_2;i Q_2)
\]
and hence
\[\text{rank image}(HF_-^*(L_1,L_2;i Q_1) \to HF_-^*(L_1,L_2;i Q_2))\]
is bounded above by $P(\lambda_i)$.
Using this fact along with the key claim,
we get $\Gamma(L_1,L_2;\K) \leq \Gamma(P)$
and hence Step 1 is complete.

\bigskip

{\bf Step 2:}
Now suppose that there is an open set $U \subset \widehat{M}$
containing $M$ and another open set $U_1 \subset U$
disjoint from $M$ so that:
\begin{itemize}
	\item $U \setminus U_1$ is compact,
	\item there is a fibration $\pi : U_1 \to V_1$
	admitting a $U(1)$ structure group so that
	the fibers are symplectomorphic to annuli
	$A_{b,B} := \{b < r < B\} \subset \C$
	and where the associated
	$U(1)$ action rotates these fibers and
	\item the symplectic vectors orthogonal to the
	fibers of $\pi$ give us an Ehresmann connection
	compatible with this $U(1)$ structure group.
\end{itemize}

We will show in this step that
for any transverse admissible Lagrangians
$L_1,L_2$ in $\widehat{M}$,
$\Gamma(L_1,L_2;\K) \leq 1$.

This is done using the main result in Step 1.
Because $\pi$ has a $U(1)$
structure group, we have a natural function
$r_1 : U_1 \to \R$, whose restriction to each
fiber $A_{b,B}$ is $r$.
Now the level sets $r_1^{-1}(c)$
of $r_1$ are coisotropic submanifolds
whose leaves are equal to the $S^1$ fibers of $\pi|_{r_1^{-1}(c)}$.
We can perturb $U$ by a generic Hamiltonian so that
$C_1 := r_1^{-1}(\frac{b+B}{2})$ is transverse to
$L_1$ and $L_2$.
One can also perturb $U$
so that $\pi|_{L_1 \cap C_1}$ and
$\pi|_{L_2 \cap C_1}$ are
smooth immersions transverse to each other
with isolated intersection points.
Let $D$ be the number of these intersection points.

We now let $H : \widehat{M} \to \R$ be a Hamiltonian
so that $C_1 = H^{-1}(c)$ for some $c$,
$H^{-1}((-\infty,c])$ is compact and so that
$H$ is
equal to $\frac{1}{2}r^2$ or $-\frac{1}{2}r^2$ near $C_1$.
Then for any level set of $H$ near $C_1$,
we can explicitly compute the intersection
points
$\phi^H_\lambda(L_1) \cap L_2$.
This is because the flow of $X_H$ near $C_1$
is just the $U(1)$ action or its inverse.
In particular, in a given level set
$H^{-1}(c')$ of $H$
for $c'$ near $c$, the number of intersection points
$H^{-1}(c') \cap \cup_{\lambda_1 \leq \lambda} (\phi^H_{\lambda_1}(L_1) \cap L_2)$
is at most $2\pi \lambda D$.
Hence we can use the main result in Step 1
to show that
$\Gamma(L_1,L_2;\K) \leq 1$.

\bigskip

{\bf Step 3:}
In this step we will construct our neighborhood
$U$ from Step 2.
We first take our cone singularity and blow
it up at the origin.
The resulting resolution is equal to a
Hermitian line bundle $L$
over a projective variety $X$
with a natural symplectic structure near
the zero section.
Because it is a Hermitian line bundle,
it has a natural radial function
$r : L \to \R$ whose level sets near the zero
section are contact hypersurfaces.
The boundary of $M$ is contactomorphic to
$r^{-1}(\epsilon)$ for some small $\epsilon>0$.
Hence we can embed the annulus bundle
$r^{-1}((2\epsilon,3\epsilon))$
symplectically into $\widehat{M}$
so that it is disjoint from $M$.
This is our associated annulus bundle $U_1$.

One can show that $\widehat{M} \setminus U_1$
has two connected components.
One of these is compact and so let $K$ be this
compact component.
Then we define $U := K \cup U_1$.
One then can apply Step 2 to show that
$\Gamma(L_1,L_2;\K) \leq 1$.

\bigskip

{\bf Comments on the more general case:}
We will now suppose that either
\begin{enumerate}[label=(\alph*)]
\item \label{item:linkcase}
$\partial M$ is contactomorphic to the link
of a general singularity
\item \label{item:smoothcase}
or $\widehat{M}$
is symplectomorphic to a smooth affine variety.
\end{enumerate}

In this case, Step 1 remains unchanged.
In Step 2, instead of having a single annulus bundle,
we get multiple annulus bundles.
The resulting $U(1)$ actions of the annulus bundles
commute with each other
(see Definition \ref{defn:nicelystratified}).
In Step 3 we construct these annulus bundles as follows:

In case \ref{item:linkcase},
we resolve the singularity at $0$.
Let $E_1,\cdots,E_l$ be the exceptional divisors.
There is an annulus bundle corresponding
to each $E_i$.
The point is that a neighborhood
of $E_i$ is a disk bundle.
We remove a smaller disk subbundle to get an
annulus bundle, which we restrict to
$E_i \setminus \cup_{j \neq i} E_i$.
Finally we symplectically embed
these annulus bundles into $\widehat{M}$.

In case \ref{item:smoothcase}
we embed our smooth affine variety into
a smooth projective variety so that the complement
is a union of smooth normal crossing divisors.
Each such divisor gives us an annulus bundle
in a similar way.

\section{Bounds for the growth rate of wrapped Floer cohomology} \label{section:growthratebounds}

Let $(M,\theta_M)$ be a Liouville domain.
Let $L_1,L_2$ be admissible Lagrangians inside $\widehat{M}$.
Let $P : [0,\infty) \to [0,\infty)$ be any function.
%
Let $L_1,L_2 \subset \widehat{M}$ be transversally intersecting admissible Lagrangians.

\begin{defn}
We say that a Hamiltonian $H : \widehat{M} \to \R$ is {\bf $(L_1,L_2,P)$-bounded} if there are some constants $C_H \in \R$,
$\delta_H >0$ so that:
\begin{PB}

\item \label{item:pbstart} $H^{-1}((-\infty,C_H \pm \delta_H])$ is compact and contains $M$.
\item  \label{item:everythingtransverse}
For all $C \in [C_H - \delta_H,C_H+\delta_H]$,
$C$ is a regular value of $H$ and
$L_1$ and $L_2$ intersect $H^{-1}(C)$ transversally.
\item \label{item:transverseintersection}
 For all $\lambda \geq 0$,
we have $\phi^H_\lambda(L_1)$ and $L_2$ intersect transversally inside a small fixed neighborhood of $H^{-1}([C_H - \delta_H,C_H+\delta_H])$
and the number of such intersection points is bounded above by $P(\lambda)$.
\item \label{item:flowlinebound} For all $C \in [C_H - \delta_H,C_H+\delta_H]$,
the number of flowlines of  $X_H$ inside $H^{-1}(C)$ of length $\leq \lambda$ starting on $L_1$ and ending on $L_2$
is bounded above by $P(\lambda)$. The length of a flowline is defined to be the time it takes to flow from start to finish.
\end{PB}
\end{defn}

Note that $H$ does not have to be admissible in the above definition.
The main theorem of this section is the following:
\begin{theorem} \label{theorem:growthratebound}
Let $P : [0,\infty) \to (0,\infty)$ be a non-decreasing function.
Let $H : \widehat{M} \to \R$ be a Hamiltonian which is $(L_1,L_2,P)$-bounded.
Then $\Gamma(L_1,L_2) \leq \Gamma(P)$.
\end{theorem}

Before we prove this theorem, we need some preliminary technical lemmas and a definition.
Good examples to keep in mind when reading these Lemmas
is when $\widehat{M} = \C$ or $\C^2$ with linear Lagrangians,
and Hamiltonians given by
a function of the radius.

\begin{lemma} \label{lemma:positiveactionhamiltonian}
Let $f : \R \to \R$ be a smooth function such that
\begin{enumerate}
\item $f,f',f'' \geq 0$ and for large $x$, $f'(x) = 1$.
\item $f(x) = 0$ if and only if $x \leq 1$.
\end{enumerate}
Let $Q_1,Q_2 : \widehat{M} \to \R$
be two Hamiltonians such that
for $j=1,2$, there are constants $C_j,\kappa_j > 0$
such that $Q_j|_{[1,\infty) \times \partial M} = \kappa_j f(r_M-C_j)$
where $r_M$ is the cylindrical coordinate,
and where $Q_j|_M = 0$.
We also assume $Q_1 \leq Q_2$.
Define $b_i := \text{rank image}(HF^*_-(L_1,L_2;i Q_1) \to HF^*_-(L_1,L_2;i Q_2))$ for each $i \in \N$.
Then: $\Gamma(L_1,L_2) \leq \text{max}(0,\Gamma((b_i)_{i \in \N}))$. 
\end{lemma}

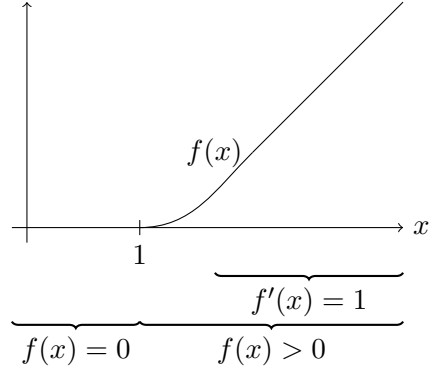
\begin{center}
\begin{tikzpicture}[domain=0:4]
   \draw[->] (-0.2,0) -- (5,0) node[right] {$x$};
   \draw[->] (0,-0.2) -- (0,3);
    \draw (1.5,0) to[out=0,in=225] (3,1) node[left] {$f(x)$};
  \draw[-] (3,1) -- (5,3);
    \draw (1.5,0.1) -- (1.5,-0.1) node[below,scale=1] {$1$};
\draw [
    thick,
    decoration={
        brace,
        mirror,
        raise=0.5cm
    },
    decorate
] (-0.2,-0.75) -- (1.5,-0.75) 
node [pos=0.5,anchor=north,yshift=-0.55cm] {$f(x)=0$};

\draw [
    thick,
    decoration={
        brace,
        mirror,
        raise=0.5cm
    },
    decorate
] (1.5,-0.75) -- (5,-0.75)
node [pos=0.5,anchor=north,yshift=-0.55cm] {$f(x)>0$}; 

\draw [
    thick,
    decoration={
        brace,
        mirror,
        raise=0.5cm
    },
    decorate
] (2.5,-0.1) -- (5,-0.1)
node [pos=0.5,anchor=north,yshift=-0.55cm] {$f'(x)=1$};

\end{tikzpicture}
\end{center}

\proof of Lemma \ref{lemma:positiveactionhamiltonian}.
Because \[Q_j -\theta_M(X_{Q_j}) = f(r_M-C_j) - (r_M-C_j)f'(r_M-C_j)\]
\[= -\int_1^{r_M} (t-C_j) f''(t-C_j) dt \leq 0\]
and because $f_{L_j} = 0$ outside $M$ where $f_{L_j}$ is the function associated to $L_j$ for $j = 1,2$ we have that
the action of every intersection point
$\phi^{Q_j}_i(L_1) \cap L_2$ outside $M \subset \widehat{M}$
is non-positive for all $i \in \N$.
Because $Q_j = 0$ inside $M$,
we get that $\phi^{Q_j}_i(L_1) = L_1$ inside $M$
and hence $\phi^{Q_j}_i(L_1)$ and $L_2$
intersect transversally inside $M$
and the number of such intersection points
is bounded above by a constant $D \geq 0$
which is independent of $i$.
Hence by property \ref{item:upperbound}
we get that
$HF^*_{(0,\infty)}(L_1,L_2,i Q_j)$ has rank bounded above by $D$ due to the fact that
intersection points outside $M$ have non-positive action.
By \ref{item:actionmorphismproperty}
this means that the kernel and cokernel of the natural action morphism
\[HF^*(L_1,L_2;i Q_j) \to HF^*_-(L_1,L_2;i Q_j)\] has rank bounded above by $D$ for all $i \geq 0$.
Because $Q_1 \leq Q_2$, we get that
$C_1 \geq C_2$ and so
 $Q_1 \leq \frac{\kappa_1}{\kappa_2} Q_2$.
The inequality $Q_1 \leq Q_2$ also tells us that
$\kappa_1 \leq \kappa_2$ and so $\frac{\kappa_1}{\kappa_2} \leq 1$.
So for all $i \in \N$, we have the following commutative diagram by property \ref{item:commutativityproperty}:
\[
\xy
(0,0)*{}="A"; (60,0)*{}="B";
(0,-20)*{}="C"; (60,-20)*{}="D";
(0,-40)*{}="E"; (60,-40)*{}="F";
"A" *{HF^*(L_1,L_2,i Q_1)};
"B" *{HF^*_-(L_1,L_2,i Q_1)};
"C" *{HF^*(L_1,L_2,i \frac{\kappa_1}{\kappa_2} Q_2)};
"D" *{HF^*_-(L_1,L_2,i \frac{\kappa_1}{\kappa_2} Q_2)};
"E" *{HF^*(L_1,L_2,i Q_2)};
"F" *{HF^*_-(L_1,L_2,i Q_2)};
{\ar@{->} "A"+(18,0)*{};"B"-(18,0)*{}};
{\ar@{->} "C"+(18,0)*{};"D"-(18,0)*{}};
{\ar@{->} "A"+(0,-4)*{};"C"+(0,4)*{}};
{\ar@{->} "B"+(0,-4)*{};"D"+(0,4)*{}};
{\ar@{->} "E"+(18,0)*{};"F"-(18,0)*{}};
{\ar@{->} "C"+(0,-4)*{};"E"+(0,4)*{}};
{\ar@{->} "D"+(0,-4)*{};"F"+(0,4)*{}};
%
"A"+(4,-10) *{a};
\endxy
\]
All the horizontal maps have kernels and cokernals of rank bounded above $D$.
Because $i Q_1$ and $i \frac{\kappa_1}{\kappa_2} Q_2$
have the same slope, we have
by Lemma \ref{lemma:hamiltonianswithsameslope}
that the map $a$ is an isomorphism.
This implies that the rank of the image of
$HF^*(L_1,L_2,i \frac{\kappa_1}{\kappa_2} Q_2) \to HF^*(L_1,L_2,i Q_2)$
is less than or equal to
\[b_i + 2D = \text{rank image}\Big(HF^*_-(L_1,L_2,i Q_1) \to HF^*_-(L_1,L_2,i Q_2)\Big)+ 2D.\]
Hence $b_i + 2D$ is bounded below by:
\[w_i := \text{rank image}\Big(HF^*(L_1,L_2,i \frac{\kappa_1}{\kappa_2} Q_2) \to HW^*(L_1,L_2)\Big).\]
Because
$w(\lambda) := \text{rank image}\Big(HF^*(L_1,L_2,\lambda \frac{\kappa_1}{\kappa_2} Q_2) \to HW^*(L_1,L_2)\Big)$
is a non-decreasing function of $\lambda \in \R$
by \ref{item:continuationmorphismproperty}, we get that
$\Gamma(L_1,L_2) = \Gamma(w(\lambda)) = \Gamma((w_i)_{i \in \N})$.
Putting everything together, we get $\Gamma(L_1,L_2) = \Gamma((w_i)_{i \in \N}) \leq \Gamma((b_i + 2D)_{i \in \N}) \leq \text{max}(0,\Gamma((b_i)_{i \in \N}))$.
\qed

\bigskip

\begin{lemma} \label{lemma:constanthamiltonianbound}
Let $H_i : \widehat{M} \to \R$ be a non-decreasing sequence of Hamiltonians and $(p_i)_{i \in \N}$
a sequence of positive reals so that:
\begin{enumerate}
\item \label{item:constantoutsidecompactset}
There is a linear function $L : \N \to \R$,
a compact set $K \subset \widehat{M}$
and a constant $c> 0$
so that
$H_i$ is a positive constant $\kappa_i$ outside $K$
with $\kappa_i \in (L(i)-c,L(i)+c)$ for all $i \in \N$.
\item \label{item:boundonintersectionptsandtransverality}
The number of intersection points $\phi^{H_i}_1(L_1) \cap L_2$ is bounded above
by $p_i$ and these are all transverse intersection points.
\item \label{item:zerooutnearM}
$H_i = 0$ on some fixed neighborhood of $M \subset \widehat{M}$.
\end{enumerate}
Then $\Gamma(L_1,L_2) \leq \Gamma((p_i)_{i \in \N})$.
\end{lemma}
\proof of Lemma \ref{lemma:constanthamiltonianbound}.
Let $f$ be the function described in the statement of Lemma \ref{lemma:positiveactionhamiltonian}.
Let $\alpha>0$ be a very small constant and $C \geq 1$ a sufficiently large constant so that:
\begin{enumerate}
\item $\{r_M \geq C\} \subset \widehat{M} \setminus K$.
\item \label{item:smallenoughham}
We have $i \alpha \left(f(r_M - C) - r_M f'(r_M - C) \right)  + \kappa_i > 0$ for all $i \in \N$ and $x \in \{r_M \geq C\}$.
This is possible because $\kappa_i > \text{max}(0,L(i)-c)$
for all $x \in \{r_M \geq C\}$.
\end{enumerate}
Define $K_i := H_i + i \alpha f(r_M -C)$.
Property (\ref{item:smallenoughham})
tells us that every intersection point $\phi^{K_i}(L_1) \cap L_2$ in the region $r_M \geq C$ has strictly positive action.

Because $L(i)$ is linear and $H_i = 0$ on some neighborhood of $M$ we can
find constants $\alpha_1,\alpha_2 > 0$ and $C_1,C_2 \geq 1$ so that:
$Q_j :=  \alpha_j f(r_M - C_j)$ satisfies:
$iQ_2 < K_i < iQ_1$ for all $i \in \N$ and $j=1,2$. 
This gives us continuation maps:
\[ HF^*_-(L_1,L_2;iQ_1) \to HF^*_-(L_1,L_2;K_i) \to HF^*_-(L_1,L_2;iQ_1)\]
by \ref{item:continuationmorphismproperty}
which implies that:
\[a_i := \text{rank image}\left(HF^*_-(L_1,L_2;iQ_1) \to HF^*_-(L_1,L_2;iQ_2)\right)\]
is bounded above by the rank of $HF^*_-(L_1,L_2;K_i)$.
Because
$K_i = H_i$ outside $\{r_i \geq C\}$
and all intersection points
$\phi_1^{K_i}(L_1) \cap L_2$ in the region $\{r_M \geq C\}$ have strictly positive action,
we have that the rank of $HF^*_-(L_1,L_2;K_i)$ is bounded above by $p_i$ by \ref{item:upperbound}.
By Lemma \ref{lemma:positiveactionhamiltonian} we then get:
\[\Gamma(L_1,L_2) \leq \text{max}(0,\Gamma((a_i)_{i \in \N}))  \leq \Gamma((p_i)_{i \in \N})).\]
\qed

\proof of Theorem \ref{theorem:growthratebound}.
We will use the Lemma \ref{lemma:constanthamiltonianbound} to give us our upper bound.
What we will do is construct a sequence of
functions $g_i : \R \to [0,\infty)$ so that the Hamiltonians
$H_i := g_i \circ H$ and the sequence $(P(\lambda_i))_{i \in \N}$ for some appropriate approximately linear sequence $(\lambda_i)_{i \in \N}$ satisfy
the conditions of Lemma \ref{lemma:constanthamiltonianbound}
and hence giving us our result.

By definition there are constants $C_H \in \R$, $\delta_H >0$ so that the conditions
\ref{item:pbstart}, \ref{item:everythingtransverse}, \ref{item:transverseintersection} and \ref{item:flowlinebound}
are satisfied.
For each $i>0$ we let $g_i : \R \to [0,\infty)$
be a smooth function (pictured below) satisfying:
\begin{gi}
\item \label{item:derivativesandconstantsatinfinity}
$g'_i \geq 0$ and $g_i \geq g_{i-1}$ for all $i \in \N$.
Also for all $i \in \N$, $g_i(x) = 0$ for $x \leq C_H - \delta_H$ and
$g_i(x) = \kappa_i$ for $x \geq C_H + \delta_H$ where
$\kappa_i \in (L(i)-c,L(i)+c)$
for some constant $c>0$
and positive linear function $L : \N \to \R$.
\item $g_i''(x) \neq 0$ if and only if $x$ is in an open set ${\mathcal N}_i$ which is a union of two open intervals whose total length is at most $\frac{\delta_H}{2i}$.
The closure of ${\mathcal N}_i$ contains
$C_H \pm \delta_H$
(in other words, the ends of these intervals
touch $C_H \pm \delta_H$).
\item $g''_i(x) \geq 0$ if $x < C_H$ and $g''_i \leq 0$ if $x > C_H$.
\item \label{item:slopegeneric}
$g'_i(C_H) = \lambda_i$ where $\lambda_i$ is a bounded distance from a positive linear function in $i$.
We also assume that $\lambda_i$ is generic enough so that there are no flowlines of $X_H$ starting on $L_1$ and finishing on $L_2$
inside $H^{-1}(C_H \pm \delta_H)$ of length $\lambda_i$.
This can be done because the set of lengths of such flowlines inside $H^{-1}(C_H \pm \delta_H)$ is discrete by property
\ref{item:flowlinebound}.
\end{gi}
The function $g_i$ will also satisfy the additional
technical property (\ref{eqn:additionalproperty})
which we state later.
We will also might need to shrink ${\mathcal N}_i$ as well later on.
But none of these changes to $g_i$ will affect \ref{item:derivativesandconstantsatinfinity}-\ref{item:slopegeneric}.
Here is a picture of $g_i$:

\begin{center}
\begin{tikzpicture}[domain=0:4]
   \draw[->] (-0.2,0) -- (9,0) node[right] {$x$};
   \draw[->] (0,-1.2) -- (0,5) node[above] {$g_i(x)$};
    \draw (2,0) to[out=0,in=236.44] (3,0.25);
    \draw (3,0.25) to (6,3.75);
    \draw (6,3.75) to[out=56.44,in=180] (7,4.2);
    \draw (7,4.2) to (9,4.2);
    \draw (2.3,0.1) -- (2.3,-0.1) node[below,scale=1] {$C_H-\delta_H$};
    \draw (4.65,0.1) -- (4.65,-0.1) node[below,scale=1] {$C_H$};
    \draw (7,0.1) -- (7,-0.1) node[below,scale=1] {$C_H+\delta_H$};
   \draw (0.1,4.2) -- (-0.1,4.2) node[left] {$\kappa_i$};
   \node at (4,3) {Slope $ = \lambda_i$};
\draw [
    thick,
    decoration={
        brace,
        mirror,
        raise=0.5cm
    },
    decorate
] (2.3,-0.75) -- (3,-0.75) 
node [pos=0.5,anchor=north,yshift=-0.55cm] {${\mathcal N}_i$}; 
\node[] at (2.55,-1) {$g''(x)>0$};
\draw [
    thick,
    decoration={
        brace,
        mirror,
        raise=0.5cm
    },
    decorate
] (6,-0.75) -- (7,-0.75) 
node [pos=0.5,anchor=north,yshift=-0.55cm] {${\mathcal N}_i$}; 
\node[] at (6.5,-1) {$g''(x)<0$};
\end{tikzpicture}
\end{center}

We define $H_i := g_i \circ H$
and hence $X_{H_i}=(g'_i \circ H)X_H$.
Condition \ref{item:pbstart} and \ref{item:derivativesandconstantsatinfinity} ensures that $H_i$ is constant outside a fixed compact set
and also zero on some open set containing $M$
and hence satisfies conditions
(\ref{item:constantoutsidecompactset})
and
(\ref{item:zerooutnearM})
of Lemma
\ref{lemma:constanthamiltonianbound}.
It remains to show that condition
(\ref{item:boundonintersectionptsandtransverality})
holds from this Lemma.

We will now show that $\phi_1^{H_i}(L_1)$ is transverse to $L_2$ for all $i$.
Outside $H^{-1}( (C_H - \delta_H, C_H + \delta_H) )$ we have that $H_i$ is constant and so $\phi_1^{H_i}(L_1)$ is transverse to $L_2$ in this region as we have assumed that $L_1$ is transverse to $L_2$.
In the region $H^{-1}(( C_H - \delta_H, C_H + \delta_H)  \setminus \mathcal{N}_i)$ we have that $X_{H_i} = \lambda_i X_H$
and so by \ref{item:transverseintersection} we have that $\phi_1^{H_i}(L_1)$ is transverse to $L_2$ in this region.

So we only need to consider $\phi^{H_i}_1(L_1)$ and $L_2$  inside
$H^{-1}({\mathcal N}_i)$.
Because $\phi^H_\lambda(L_1)$
is transverse to $L_2$
in $H^{-1}({\mathcal N}_i)$
for all $\lambda \geq 0$ by \ref{item:transverseintersection},
we get that the manifolds
$\left\{ (\lambda,x) \in (0,\lambda_i] \times
H^{-1}({\mathcal N}_i)
\left| x \in \phi^H_\lambda(L_1) \right.\right\}$
and
$(0,\lambda_i] \times L_2$
intersect transversally inside $(0,\lambda_i] \times H^{-1}({\mathcal N}_i)$.
Hence 
\[\Lambda = \left\{ (\lambda,x) \in (0,\lambda_i] \times
H^{-1}({\mathcal N}_i)
\left| x \in \phi^H_\lambda(L_1) \cap L_2 \right.\right\}
\]
is a proper $1$-dimensional submanifold of
$(0,\lambda_i] \times H^{-1}({\mathcal N}_i)$.
Let $p_1 : \Lambda \to (0,\lambda_i]$,
$p_2 : \Lambda \to H^{-1}({\mathcal N}_i)$ be 
the natural projection maps to $(0,\lambda_i]$ and
$H^{-1}({\mathcal N}_i)$ respectively.
Let $h:=p_2^*(H|_{H^{-1}({\mathcal N}_i)})$.

We can make sure that
${\mathcal N}_i$ is small enough for each $i$
so that
$M_i := \sup(p_1(\Lambda)) < \lambda_i$
by property \ref{item:slopegeneric}
combined with the fact that
the set of lengths of flowlines of $X_H$
inside $H^{-1}(C_H \pm \delta_H)$
is discrete by \ref{item:flowlinebound}
and $m_i := \inf(p_1(\Lambda)) > 0$
by the last part of \ref{item:everythingtransverse}.

Choose an open subset
${\mathcal N}^*_i \subset {\mathcal N}_i$
which is also a union of two non-empty intervals,
one in each connected component of ${\mathcal N}_i$,
such that each $x \in H^{-1}({\mathcal N}^*_i)$
is a regular value of $h$
(i.e. $dh \neq 0$ along $h^{-1}(x)$).

We can modify $g_i$ so that it satisfies the following
additional property:
\begin{equation} \label{eqn:additionalproperty}
\text{if} \quad g'(x) \in [m_i,M_i]
\ \text{then} \ x \in {\mathcal N}^*_i
\ \text{and} \ |g_i''(x)| > \text{max}\left| \left. \frac{dp_1}{dh}\right|_{h^{-1}(x)} \right|.
\tag{gi5}
\end{equation}
Here is a graph of $g'_i$:
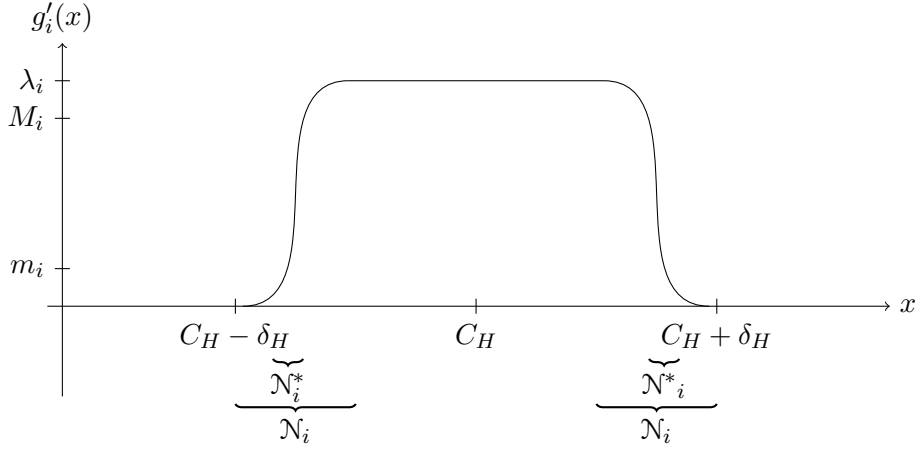
\begin{center}
\begin{tikzpicture}[domain=0:4]
   \draw[->] (-0.2,0) -- (11,0) node[right] {$x$};
   \draw[->] (0,-1.2) -- (0,3.5) node[above] {$g'_i(x)$};
    \draw (2.4,0) to[out=0,in=180] (3.8,3);
    \draw (3.8,3) to (7.2,3);
    \draw (7.2,3) to[out=0,in=180] (8.6,0);
    \draw (2.3,0.1) -- (2.3,-0.1) node[below,scale=1] {$C_H-\delta_H$};
    \draw (5.5,0.1) -- (5.5,-0.1) node[below,scale=1] {$C_H$};
    \draw (8.7,0.1) -- (8.7,-0.1) node[below,scale=1] {$C_H+\delta_H$};
   \draw (0.1,3) -- (-0.1,3) node[left] {$\lambda_i$};
   \draw (0.1,2.5) -- (-0.1,2.5) node[left] {$M_i$};
   \draw (0.1,0.5) -- (-0.1,0.5) node[left] {$m_i$};
\draw [
    thick,
    decoration={
        brace,
        mirror,
        raise=0.5cm
    },
    decorate
] (2.3,-0.8) -- (3.9,-0.8) 
node [pos=0.5,anchor=north,yshift=-0.55cm] {${\mathcal N}_i$}; 
\draw [
    thick,
    decoration={
        brace,
        mirror,
        raise=0.5cm
    },
    decorate
] (2.8,-0.2) -- (3.2,-0.2) 
node [pos=0.5,anchor=north,yshift=-0.55cm] {${\mathcal N}^*_i$}; 
\draw [
    thick,
    decoration={
        brace,
        mirror,
        raise=0.5cm
    },
    decorate
] (7.1,-0.8) -- (8.7,-0.8) 
node [pos=0.5,anchor=north,yshift=-0.55cm] {${\mathcal N}_i$}; 
\draw [
    thick,
    decoration={
        brace,
        mirror,
        raise=0.5cm
    },
    decorate
] (7.8,-0.2) -- (8.2,-0.2) 
node [pos=0.5,anchor=north,yshift=-0.55cm] {${\mathcal N^*}_i$}; 
\end{tikzpicture}
\end{center}

Let $V$ be a non-zero vector tangent to
$L_1$ at a point $p \in H^{-1}({\mathcal N}_i)$.
We wish to show that
$\check{V} := D\phi^{H_i}_1(V) \notin TL_2$ for all $i \in \N$.
We have two cases to consider:

{\it Case 1}: $V$ is not tangent to any level set of $H$.

{\it Case 2}: $V$ is tangent to a level set of $H$.

{\it Case 1}:
Suppose (for a contradiction)
that $\check{V} \in TL_2$.
After rescaling $V$ we can assume that
$dH(V) = 1$ and hence $dH(\check{V})=1$. Because the length of any flowline
of $X_H$ starting at $p$ and ending in $L_2$
is in $[m_i,M_i]$, we get that $g'(H(p)) \in [m_i,M_i]$
and hence $p \in H^{-1}({\mathcal N}_i^*)$.
Therefore $\check{V}$ is in the image of
$Dp_2|_{p_2^{-1}({\mathcal N}_i^*)}$.
Because $p_2|_{p_2^{-1}({\mathcal N}_i^*)}$ is an
immersion, there is a
unique vector $\widetilde{V} \in T\Lambda$
so that $Dp_2(\widetilde{V}) = \check{V}$.
Because $dH(\check{V})=1$, we get
$dp_1(\widetilde{V}) = \frac{dp_1}{dh}$.
Because the flowlines of $H_i$ of length $1$
are equal, up to reparameterization,
to the flowlines of $H$ of length $g'(H)$,
we get that $g''(H(p)) =dp_1(\widetilde{V})$.
But this contradicts (\ref{eqn:additionalproperty}).

{\it Case 2}: 
We now consider the case when $V$ is tangent to
$L_1$ and also tangent to a level set $H^{-1}(C)$
of $H$ inside
$H^{-1}({\mathcal N}_i)$.
Inside $H^{-1}(C)$
we have $X_{H_i}=g'_i(C)X_H$ and so by property
\ref{item:everythingtransverse},
$\phi^{H_i}_1(L_1) \cap H^{-1}(C)$ and $L_2 \cap H^{-1}(C)$ are submanifolds of $H^{-1}(C)$.
Also $X_{H_i}$ is a constant multiple of
$X_H$ inside $H^{-1}(C)$.
Hence by
\ref{item:transverseintersection}
we get that at each intersection point
$q \in \phi^{H_i}_1(L_1) \cap L_2 \cap H^{-1}(C)$,
the tangent spaces at $q$ of
$\phi^{H_i}_1(L_1) \cap H^{-1}(C)$ and $L_2 \cap H^{-1}(C)$ intersect in one point.
Hence $\check{V} \notin TL_2$.
Putting everything together we get that $\phi_1^{H_i}(L_1)$ is transverse to $L_2$ for all $i$.

We now need a bound on the number of intersection points
$\phi^{H_i}_1(L_1) \cap L_2$.
The number of intersection points $\phi^{H_i}_1(L_1) \cap L_2$ in the complement of $H^{-1}({\mathcal N}_i)$
is bounded above by $P(\lambda_i) + D$
where $D$ is a constant by \ref{item:transverseintersection} because $X_{H_i}$ is
equal to $\lambda_i X_H$ inside $H^{-1}( [C_H - \delta_H, C_H + \delta_H ]) \setminus {\mathcal N}_i$ and zero outside the union of this region and ${\mathcal N}_i$.
Choose two points $a_-,a_+ \in {\mathcal N}^*_i$,
one in each connected component.
Property (\ref{eqn:additionalproperty})
combined with the fact that
\begin{itemize}
	\item $\Lambda$ is a proper submanifold of $(0,\lambda_i] \times H^{-1}({\mathcal N}_i)$,
	\item $\text{image}(p_1) \subset [m_i,M_i] \subset (0,\lambda_i)$
\end{itemize}
implies that
there is a one to one correspondence between length $1$ flowlines of $X_{H_i}$
from $L_1$ to $L_2$ inside $H^{-1}({\mathcal N}_i)$ and flowlines of $X_H$ from $L_1$ to $L_2$ of length at most $\lambda_i$ inside $H^{-1}(a_+) \cup H^{-1}(a_-)$.
Hence the number of length $1$ flowlines of $X_{H_i}$ in $H^{-1}({\mathcal N}_i)$
is equal to the number of flowlines of $X_H$ inside $H^{-1}(a_-)$ plus the number flowlines of $X_H$
in $H^{-1}(a_+)$ all of length $\leq \lambda_i$.
Hence there are at most $2P(\lambda_i)$ length $1$ flowlines of $X_{H_i}$ in $H^{-1}({\mathcal N}_i)$ by \ref{item:flowlinebound}.
So there are at most $2P(\lambda_i)(P(\lambda_i) + D)$ length $1$ flowlines of $X_{H_i}$ starting at $L_1$ and ending at $L_2$.

Using this fact combined with the fact that
$\phi^{H_i}_1(L_1)$ and $L_2$ intersect transversally
for all $i$ we get that condition
(\ref{item:boundonintersectionptsandtransverality})
holds in Lemma \ref{lemma:constanthamiltonianbound}
for $H_i$ and the sequence $p_i := P(\lambda_i)$.
Hence by Lemma \ref{lemma:constanthamiltonianbound},
$\Gamma(L_1,L_2) \leq \Gamma(((2P(\lambda_i)(P(\lambda_i) + D))_{i \in \N})) = \Gamma((P(\lambda_i))_{i \in \N}) = \Gamma(P)$.
\qed

\section{Growth rates and Compatible Annulus Bundles} \label{section:nicelystratifiedsection}

Let $(r,\vartheta)$ be the standard polar coordinates in $\C$ and let $A_{b_1,b_2}$
be the open annulus equal to $\{ b_1 < r < b_2\}$ with the standard symplectic structure $\frac{1}{2} d(r^2) \wedge d\vartheta$.

\begin{defn} \label{defn:nicelystratified} \label{definition:nicelystratifiedatinfinity}
Let $(M,\theta_M)$ be a Liouville domain. We say that $\widehat{M}$  {\bf admits compatible annulus bundles at infinity of codimension $C_D \in \N$}
if 
for each $I \subset \{1,\cdots,l\}$, there are open subsets $U_I \subset \widehat{M}$,
 manifolds $V_I$ and smooth fibrations $\pi_I : U_I \twoheadrightarrow V_I$ satisfying the following properties:
\begin{NS}
\item \label{item:containsa}
We have $U_{I \cup J} = U_I \cap U_J$ for all
$I,J \subset \{1,\cdots,l\}$.
Also $U_\emptyset \setminus \cup_i U_{\{i\}}$ contains $M$.
If $|I| > C_D$ then $U_I = \emptyset$.
\item \label{item:structuregroupcondition}
Each fiber of $\pi_I$
is a symplectic submanifold of $U_I$ for each $I$.
If $I = \emptyset$ then
$\pi_\emptyset : U_{\emptyset} \to V_{\emptyset}$ is a diffeomorphism.
For non-empty $I$, the fiber bundle $\pi_I$ has a $U(1)^{|I|}$ structure group where the fibers are symplectomorphic to a product of symplectic manifolds
$\prod_{j \in I} A_j$
and where:
\begin{itemize}
\item $A_j = A_{b,B}$ for some fixed constants $0 < b < B$ (independent of $j$)
\item and the structure group $U(1)^{|I|}$ acts on $\prod_{j \in I} A_j$ by $(e^{i t_j})_{j \in I} \cdot (x_j)_{j \in I} = (e^{i t_j}x_j)_{j \in I}$.
\end{itemize}

The set of vectors in $U_I$ which are symplectically orthogonal to the fibers of $\pi_I$
give us an Ehresmann connection respecting the above $U(1)^{|I|}$ structure group.

\item \label{item:fiberproperty}
For $J \subset I$, the fibers of $\pi_J|_{U_I}$ are contained in the fibers of $\pi_I$ and in any trivialization $\prod_{j \in I} A_j$
of any fiber of $\pi_I$, the fibers of $\pi_J|_{U_I}$ are of the form
$\prod_{j \in J} A_j \times \prod_{j \in I \setminus J} \{z_j\}$ for points $z_j \in A_j$.
Also for all $I \subset \{1,\cdots,l\}$ we have $\pi_I(U_I \setminus \cup_{i \notin I} U_{\{i\}})$ is a compact subset of $V_I$.

\item \label{item:closureproperty}
We define $r_i : U_{\{i\}} \to (b,B)$ to be the unique function
whose restriction to each fiber of $\pi_{\{i\}}$
is equal to the radial coordinate in $A_i = A_{b,B} \subset \C$
in any $U(1)$ trivialization.
We call $r_i$ the $i$th {\bf radial coordinate}.
We require that the set
$\{r_i \geq r\}$
is closed in $U_\emptyset$ for each $i \in \{1,\cdots,l\}$ and
$r \in (b,B)$.
\end{NS}

The data $U_I,\pi_I,V_I$ is called a
{\bf set of compatible annulus bundles on} $\widehat{M}$.
\end{defn}

Here is a schematic picture of a particular example of compatible annulus bundles at infinity of codimension $C_D = 2$:
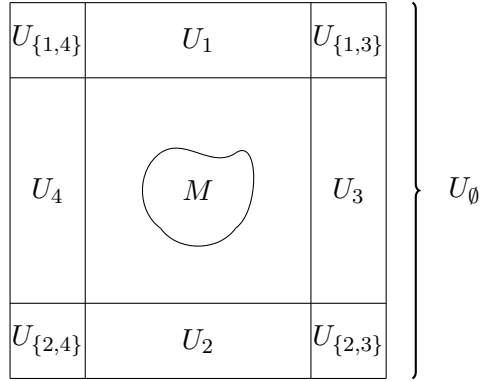
\begin{center}
\begin{tikzpicture}[domain=0:4]
\draw (0,0) -- (0,5);
\draw (0,5) -- (5,5);
\draw (5,5) -- (5,0);
\draw (5,0) -- (0,0);
\draw (1,0) -- (1,5);
\draw (0,4) -- (5,4);
\draw (4,5) -- (4,0);
\draw (5,1) -- (0,1);
\draw (2,2) to[out=305,in=235] (3,2);
\draw (3,2) to[out=35,in=35] (3,3);
\draw (3,3) to[out=215,in=35] (2,3);
\draw (2,3) to[out=215,in=145] (2,2);
\node[] at (2.5,2.5) {$M$};
\node[] at (2.5,4.5) {$U_1$};
\node[] at (2.5,0.5) {$U_2$};
\node[] at (4.5,2.5) {$U_3$};
\node[] at (0.5,2.5) {$U_4$};
\node[] at (4.5,4.5) {$U_{\{1,3\}}$};
\node[] at (4.5,0.5) {$U_{\{2,3\}}$};
\node[] at (0.5,0.5) {$U_{\{2,4\}}$};
\node[] at (0.5,4.5) {$U_{\{1,4\}}$};
\draw [
    thick,
    decoration={
        brace,
        mirror,
        raise=0.5cm
    },
    decorate
] (4.85,0) -- (4.85,5) 
node [pos=0.5,anchor=east,xshift=1.55cm] {$U_\emptyset$}; 
%
%
\end{tikzpicture}
\end{center}

{\bf Notational convention:}
If $\tau : \prod_{j \in I} A_j \to \pi_I^{-1}(q)$
is a $U(1)^{|I|}$ trivialization
of some fiber of $\pi_I$ over $q$
then we will just write:
$\prod_{j \in I} A_j$ for such a trivialization.
Each annulus $A_j$ has natural polar coordinates
$(r_j,\vartheta_j)$.
The coordinates $(r_i,\vartheta_i)_{i \in I}$
inside this trivialization will be called the
{\bf associated polar coordinates}.
Note that by abuse of notation,
$r_i$ is also the $i$th radial coordinate on $U_i$.

From now on we will assume that $\widehat{M}$ admits compatible annulus bundles at infinity of codimension $C_D$ where $U_I,V_I,\pi_I,b,B,r_i$
are from Definition \ref{definition:nicelystratifiedatinfinity} above.
Let $L_1$ and $L_2$ be two admissible Lagrangians inside $\widehat{M}$ that intersect each other transversally.
The main aim of this section is to prove:
\begin{theorem} \label{theorem:nicelystratifiedgrowth}
For every field $\K$, $\Gamma(L_1,L_2;\K) \leq C_D$.
\end{theorem}

We will write $\pi_i,U_i,V_i$ instead of $\pi_{\{i\}},U_{\{i\}},V_{\{i\}}$ respectively.
By abuse of notation
we define the set $\{r_i \leq a\} \subset U_\emptyset$ as the  union of $U_\emptyset \setminus U_i$
and $\{r_i \leq a\} \subset U_i$.

\begin{defn}
Let $p \in U_I$ and let $F_p$
be the fiber of $\pi_I$ passing through $p$.
A Lagrangian $L \subset \widehat{M}$
is called {\bf fiber radial near} $p$ if
there is a neighborhood ${\mathcal N}_p$ of $p$
and a trivialization $\prod_{j \in I} A_j$
of $F_p$
so that $L \cap {\mathcal N}_p \cap F_p = R \cap {\mathcal N}_p$
where $R$ is the real part of $\prod_{j \in I} A_j$
(viewed as an open subset of $\C^{|I|}$).

An alternative definition of this is as follows:
$L$ is fiber radial near $p$ if
there is a neighborhood ${\mathcal N}_p$ of $p$ so that
for any trivialization of $F_p$
with associated polar coordinates $(r_j,\vartheta_j)_{j \in I}$,
there are constants $(c_j)_{j \in I}$ so that
\[L \cap {\mathcal N}_p \cap F_p =
{\mathcal N}_p \cap F_p \cap \cap_{j \in I} \{ \vartheta_j = c_j \}.\]
\end{defn}

Fix some small $\epsilon_g > 0$.
We can assume that $\epsilon_g$ is generic
enough so that $L_1$ and $L_2$
are transverse to $\cap_{j \in I} \{r_j = b+\epsilon_g\}$ for each $I \subset \{1,\cdots,l\}$.
Let $C_I := \cap_{j \in I} \{r_j = b+\epsilon_g\} \setminus \cup_{j \notin I} \{r_j \geq b+\epsilon_g\}$
for $I \neq \emptyset$ and $C_\emptyset = U_\emptyset \setminus \cup_i \{r_i \geq b+\epsilon_g\}$.
These are coisotropic submanifolds such that the leaves of their associated coisotropic foliations
are given by fibers of $\pi_I|_{C_I}$.
Here $\pi_I|_{C_I}$ is a principal $U(1)^{|I|}$
bundle where the $U(1)^{|I|}$ structure group
is inherited from $\pi_I$.
The submanifolds $C_I$ are not usually properly embedded in $\widehat{M}$
(their closures are submanifolds with corners).

\begin{lemma} \label{lemma:lagrangianpertubation}
There is a $C^0$ small compactly supported Hamiltonian diffeomorphism
$\phi : \widehat{M} \to \widehat{M}$
so that for each $I \subset \{1,\cdots,k\}$,
\begin{enumerate}
\item \label{item:lemma441}
$\phi(L_1),\phi(L_2)$ are transverse to $C_I$ and
the Lagrangian immersions $\iota^\phi_{1,I} := \pi_I|_{\phi(L_1) \cap C_I}$ and $\iota^\phi_{2,I} := \pi_I|_{\phi(L_2) \cap C_I}$ are transverse to each other
and also the intersection points between $\iota^\phi_{1,I}$
and $\iota^\phi_{2,I}$ are isolated if $|I|<n$,
\item \label{item:lemma442}
and if
$x \in \text{Image}(\iota^\phi_{1,I}) \cap \text{Image}(\iota^\phi_{2,I})$
then $\phi(L_1)$ (resp. $\phi(L_2)$) is fiber radial near each point of
$\phi(L_1) \cap \pi_I^{-1}(x) \cap C_I$
(resp. $\phi(L_2) \cap \pi_I^{-1}(x) \cap C_I$).
\end{enumerate}
\end{lemma}

Before we prove this we need a definition and some preliminary lemmas
(Lemmas \ref{lemma:exactisotopyimplieshamisotopy},
\ref{lemma:transverseimmersions},
\ref{lemma:perturbinglagrangian},
\ref{lemma:patchingfunctions},
\ref{lemma:darbouxchart} and \ref{lemma:lagrangiansinstandardform}).
The definition and
Lemmas \ref{lemma:exactisotopyimplieshamisotopy},
\ref{lemma:transverseimmersions},
\ref{lemma:perturbinglagrangian}
deal with Lemma \ref{lemma:lagrangianpertubation}
part (\ref{item:lemma441})
and the other lemmas deal with part
(\ref{item:lemma442}).

\begin{defn} \label{defn:exactisotopyofimmersions}
Let $\iota : L \to S$ be any smooth map into a symplectic manifold $(S,\omega_S)$.
Then an {\bf exact isotopy of $\iota$} consists of a smooth family of
maps $\iota_t : L \to S, t \in [0,1]$ where $\iota_0 = \iota$
and $(\iota_t)^* (i_{\frac{d}{dt}(\iota_t)} \omega_S) = df_t$ for some smooth family of functions
$f_t : L \to \R$.
We call $(f_t)_{t \in [0,1]}$ {\bf functions associated to $\iota_t$}.
The {\bf support} $\iota_t$ is the set
of points $s \in S$
where $\iota_t(y) = s$
and $\frac{d}{dt}(\iota_t(y)) \neq 0$
for some $y \in L$.

%
\end{defn}

\begin{lemma} \label{lemma:exactisotopyimplieshamisotopy}
If $\iota_t : L \to S$ is an exact isotopy of $\iota$ where $(S,\omega_S)$
is a symplectic manifold and $\iota_t$ are smooth embeddings,
then $\iota_t$ is a Hamiltonian isotopy of $\iota$.
In other words, there is a Hamiltonian $H_t$
so that $\phi^{H_t}_t \circ \iota_0 = \iota_t$.
If $\frac{d}{dt}(\iota_t)$ is $C^\infty$ small,
we can ensure that $H_t$ is $C^\infty$ small.
\end{lemma}
\proof
Let $(f_t)_{t \in [0,1]}$ be functions associated to $\iota_t$.
Choose a smooth family of Hamiltonians
$H_t : S \to \R$ parameterized by $t \in \R$
so that $\iota_t^* H_t = f_t$ and so that
$dH_t = i_{\frac{d}{dt}(\iota_t)} \omega_S$ along
$\text{Image}(\iota_t)$ for all $t$.
Then this smooth family of Hamiltonians generates our exact isotopy and hence is a Hamiltonian isotopy.

If $\frac{d}{dt}(\iota_t)$ is $C^\infty$ small,
then by subtracting an appropriate smooth family of constants $c_t$ from $f_t$, we can ensure that $f_t$ is $C^\infty$ small.
Also $i_{\frac{d}{dt}(\iota_t)} \omega_S$ is $C^\infty$ small and hence
one can choose $H_t$ to be $C^\infty$ small.
\qed

\begin{lemma} \label{lemma:transverseimmersions} \label{lemma:transverselagrangianisotopies}
	Let $S$ be a symplectic manifold
	of dimension $\geq 2$, $K \subset compact$,
	and let $\iota : L \to S$, $\check{\iota} : \check{L} \to S$
	be proper Lagrangian immersions.
	Then there are $C^\infty$ small exact
	isotopies $\iota_t$ and $\check{\iota}_t$
	of $\iota$ and $\check{\iota}$ respectively
	supported near $K$ so that
	$\iota_1$ is transverse to $\check{\iota}_1$ with isolated intersection points.
	These intersection points can avoid
	any fixed finite union of hypersurfaces.
\end{lemma}
\proof
	Let $n = \frac{1}{2}\text{dim}(S)$.
	For each $x \in \text{image}(L) \cap K$
	choose a smooth family of compactly supported Hamiltonians
	$(H^x_s : S \to \R)_{s \in S^{2n-1}}$
	where $S^{2n-1}$ is the unit sphere in $\R^{2n}$
	so that the map:\[S^{2n-1} \times \R \to T_x S, \quad (s,t) \to tX_{H^x_s}|_{T_x S}\] is surjective.
	Because $\text{image}(L) \cap K$ is compact, there is a finite set of points
	$x_1,\cdots,x_l \in \text{image}(L)$ and a small relatively compact neighborhood $U$ of $K$ so that
	for all $x \in U$, there is some $i$ so that
	the map: $S^{2n-1} \times \R \to T_x S, \quad (s,t) \to tX_{H^{x_i}_s}|_{T_x S}$ is surjective.
	
	Now define a smooth family of autonomous Hamiltonians
	$(H_s : S \to \R)_{s = (s_1,\cdots,s_l) \in (S^{2n-1})^l}, \quad
	H_s = \sum_{i = 1}^l H^{x_i}_{s_i}$.
	Let $\phi^s_t$ be the time $t$ flow of $H_s$.
	Define $\Phi : L \times (S^{2n-1})^l \times \R \to S$
	by $\Phi(y,s,\tau) := \phi^s_\tau(\iota(y))$.
	Then $\Phi$ is a submersion for all $|\tau| < \epsilon$ where $\epsilon >0$ is sufficiently small.
	Let $\check{\Phi}$ be the restriction of $\Phi$ to $L \times (S^{2n})^l \times (-\epsilon,\epsilon)$.
	By Sard-Smale we can choose generic $(s,\tau)$ sufficiently close to $(0,0)$
	so that $L \times \{s\} \times \{\tau\}$ is transverse to the manifold
	$\check{\Phi}^{-1}(\text{image}(\iota|_V) \cap U)$ for all open $V \subset L$ such that $\iota|_V$ is an embedding.
	Then $\widetilde{\iota}_t : L \to S$
	defined by $\widetilde{\iota}_t(y) := \check{\Phi}(y,s,t\tau)$ is a $C^\infty$ small exact isotopy so that
	so that $\widetilde{\iota}_1$ is transverse to $\check{\iota}$ inside $U$.
	Because there are only finitely many points $(y,\check{y}) \in L \times \check{L}$ with
	$\widetilde{\iota}_1(y) = \check{\iota}(\check{y}) \in U$, it is fairly easy to find exact isotopies of $\widetilde{\iota}_1$ and $\check{\iota}$
	supported near each such $y$ and $\check{y}$ so that our intersection points become isolated
	and disjoint from any fixed finite union of hypersurfaces.
\qed

\begin{lemma} \label{lemma:perturbinglagrangian}
Let $(S,\omega_S)$ be a symplectic manifold and $C \subset S$ a
coisotropic submanifold which is not Lagrangian with the property that the canonical foliation
on $C$ is a fibration $\pi : C \twoheadrightarrow B$ with compact fibers.
Let $L,\check{L} \subset S$ be Lagrangians transverse to $C$. We do not assume that $S,C,B,L,\check{L}$ are compact.
Let $K \subset B$ be a compact subset of $B$.
Then we can perturb $L$ and $\check{L}$ by a $C^\infty$ small Hamiltonian isotopy supported inside a small neighborhood of $\pi^{-1}(K)$
so that the immersions $\pi|_{L \cap C}$ and $\pi|_{\check{L} \cap C}$ become transverse in a small neighborhood
of $K$. We can also ensure that their intersection points are isolated
in this small neighborhood
and that they avoid any fixed finite union of
hypersurfaces in $B$.
\end{lemma}
\proof of Lemma \ref{lemma:perturbinglagrangian}.
Note that $B$ has a natural symplectic structure $\omega_B$ defined as follows: if $X,Y$ are two vectors in $B$
then $\omega_B(X,Y)$ is defined as $\omega_S(\widetilde{X},\widetilde{Y})$ where $\widetilde{X},\widetilde{Y}$ are any vectors
tangent to $C$ satisfying $\pi_*(\widetilde{X})=X$ and $\pi_*(\widetilde{Y})=Y$.
We have that $\pi|_{L \cap C}, \pi|_{\check{L} \cap C}$ are Lagrangian immersions because their lifts $L \cap C, \check{L} \cap C$
are isotropic inside $S$.
By Lemma \ref{lemma:transverseimmersions},
one can find exact isotopies $\iota_t : L \cap C \to B$,
$\check{\iota}_t : \check{L} \cap C \to B$
of the immersions $\pi|_{L \cap C}$ and $\pi|_{\check{L} \cap C}$
so that $\iota_1$ and $\check{\iota}_1$ are transverse near $K$
and $\iota_t$, $\check{\iota}_t$ have support in some small neighborhood of $K$.
We can also ensure that intersection points between $\iota_1$ and $\check{\iota}_1$
 near $K$ are isolated and avoid any fixed finite union of hypersurfaces.
Choose an Ehresmann connection on $\pi$.
We define $(\widetilde{\iota}_t : L \cap C \to S)_{t \in [0,1]}$
so that $\widetilde{\iota}_t(y) \in C$
is equal to $y$ parallel transported
along the path $[0,t] \to B, \quad s \to \iota_s(y)$ with respect to this connection.
Define
$(\widetilde{\check{\iota}}_t : L \cap C \to S)_{t \in [0,1]}$
in the same way with respect to $\check{\iota}_t$.
Then
$(\widetilde{\iota}_t : L \cap C \to S)_{t \in [0,1]}$
and
$(\widetilde{\check{\iota}}_t : L \cap C \to S)_{t \in [0,1]}$
are smooth families of embeddings
so that
\begin{itemize}
\item the image of $\iota_t$ and $\check{\iota}_t$
is contained in $C$ and $\iota_t = \pi \circ \widetilde{\iota}_t$
and $\check{\iota}_t = 
\pi \circ \widetilde{{\check{\iota}}}_t$,
\item $\frac{d}{dt}(\widetilde{\iota}_t)$
and $\frac{d}{dt}(\widetilde{\check{\iota}}_t)$
are $C^\infty$ small, $\iota_0=\text{id}$
and  $\check{\iota}_0=\text{id}$,
\item $(\widetilde{\iota}_t : L \cap C \to S)_{t \in [0,1]}$
and
$(\widetilde{\check{\iota}}_t : L \cap C \to S)_{t \in [0,1]}$
are exact isotopies.

\end{itemize}
Hence by Lemma \ref{lemma:exactisotopyimplieshamisotopy},
there are $C^\infty$ small Hamiltonians
$H : [0,1] \times S \to S$ and $\check{H} : [0,1] \times S \to S$ supported in a small neighborhood of $\pi^{-1}(K)$
whose flow generates the isotopies $\widetilde{\iota}_t$,
$\widetilde{\check{\iota}}_t$ respectively.
Let $L^1 := \phi^H_1(L)$ and ${\check{L}}^1 := \phi^{\check{H}}_1(\check{L})$.
Then the immersions $\pi|_{L^1 \cap C}$ and $\pi|_{{\check{L}}^1 \cap C}$
have isolated transverse intersection points near $K$ that avoid our fixed union of hypersurfaces.
Hence $L^1$ and ${\check{L}}^1$ are our perturbations of $L$ and $\check{L}$ respectively.
\qed

\begin{lemma} \label{lemma:patchingfunctions}
Let $S_1, S_2 \subset \R^m$ be two submanifolds which intersect transversely at $0$
and let $f_1,f_2 : \R^m \to \R$ be two smooth functions so that $f_1(0) = f_2(0)$ and $df_1 = df_2$ at $0$.
Then there is a smooth function $f : {\mathcal N} \to \R$ where ${\mathcal N}$ is an open neighborhood of $0$
so that $f|_{S_i \cap {\mathcal N}} = f_i|_{S_i \cap {\mathcal N}}$ for $i=1,2$
and so that the restriction of $df$ to $S_2 \cap {\mathcal N}$ is equal to the restriction of $df_2$ to $S_2 \cap {\mathcal N}$.
\end{lemma}
\proof
This is done by looking at the graphs of $f_1$ and $f_2$ respectively and patching them together
using the exponential map with respect to a metric which makes the graph of $f_1$ into a totally geodesic submanifold.
Here are the details of the argument.

Let $\Gamma_1,\Gamma_2 \subset \R^m \times \R$ be the graphs of $f_1$ and $f_2$ respectively.
Let $\widetilde{S}_i := (S_i \times \R) \cap \Gamma_i$ for $i = 1,2$.
Choose a metric on $\R^m \times \R$ making $\widetilde{S}_1$ totally geodesic near $(0,f_1(0))$ with the additional property that it is complete
and so that $T\widetilde{S}_1$ and $T\widetilde{S}_2$ are orthogonal at $(0,f_1(0))$.
Here is how to construct such a metric:
First choose a complete metric making $\widetilde{S}_1$ totally geodesic near $0$
and then pull such a metric back by a compactly supported diffeomorphism fixing $(0,f_1(0))$ and sending
$\widetilde{S}_1$ to itself and whose linearization
at $(0,f_1(0))$ sends $T_{(0,f_1(0))} \widetilde{S}_2$ to vectors orthogonal to $\widetilde{S}_1$.
Here $T_{(0,f_1(0))} \widetilde{S}_2$ is the tangent space to $\widetilde{S}_2$ at $(0,f_1(0))$.

Let $N\widetilde{S}_2$ be the normal bundle of $\widetilde{S}_2$ in $\R^m \times \R$
(I.e. the set of vectors in $T(\R^m \times \R)|_{\widetilde{S}_2}$ which are orthogonal to $\widetilde{S}_2$
with respect to our chosen metric)
and let
$N^\delta \widetilde{S}_2 \subset N\widetilde{S}_2$ be the subset consisting of normal vectors of length $< \delta$.
Let $\exp^\delta : N^\delta S_2 \to \R^m \times \R$ be the exponential map.
Define $\widehat{N}^\delta\widetilde{S}_2 \subset N^\delta \widetilde{S}_2$
to be the subset consisting of vectors tangent to $\Gamma_2$.
For $\delta$ small enough, we have that the image of
$\exp^\delta|_{\widehat{N}^\delta\widetilde{S}_2}$ is the graph of some function
$f$ on an open set ${\mathcal N} \subset \R^m$ containing $0$.
We will also assume that ${\mathcal N}$ is small enough so that
it does not intersect
$(S_1 \cap S_2) \setminus \{0\}$.

Because
\begin{itemize}
	\item $\widetilde{S}_1$ is totally geodesic,
	\item $T_{(0,f_1(0))} \Gamma_1 = T_{(0,f_1(0))} \Gamma_2$,
	\item $\widetilde{S}_1$ and $\widetilde{S}_2$ are orthogonal at $(0,f_1(0))$,	
\end{itemize} 
we get that $f|_{S_1 \cap {\mathcal N}} = f_1|_{S_1 \cap {\mathcal N}}$.
Also because $\Gamma_2$ and the graph of $f$ have identical tangent spaces inside
 $T(\R^m \times \R)|_{\widetilde{S}_2 \cap ({\mathcal N} \times \R)}$,
we get that $f|_{S_2 \cap {\mathcal N}} = f_2|_{S_2 \cap {\mathcal N}}$
and $df = df_2$ along $S_2 \cap {\mathcal N}$.
\qed

\begin{lemma} \label{lemma:darbouxchart}
	Let $(S,\omega_S)$ be a symplectic manifold
	and let $f_1,\cdots,f_s$ be Poisson commuting
	functions so that $df_1,\cdots,df_s$
	are linearly independent at each point of $S$.
	Then for each $p \in \cap_{i=1}^s f_i^{-1}(0)$,
	there is a Darboux chart centered at $p$
	so that some of the Darboux coordinates
	are $f_1,\cdots,f_s$.
\end{lemma}
\proof
First of all we find functions
$f_{s+1},\cdots,f_n$ which Poisson
commute with $f_1,\cdots,f_s$
and each other by induction
and so that $df_1,\cdots,df_n$ are linearly independent at each point near $p$.

Suppose (by induction) we have constructed $f_1,\cdots,f_l$ for some $s \leq l < n$.
Because ${\mathcal L}_{X_{f_i}}(\omega_S) = 0$ and ${\mathcal L}_{X_{f_i}}df_j = 0$ for all $i,j$,
we get that $X_{f_1},\cdots,X_{f_l}$ commute
as vector fields.
Hence by the Frobenius integrability theorem,
there is a neighborhood $U$ of $p$
and a smooth map $P : U \to V$
whose fibers are regular of dimension $l<n$
and are tangent to $X_{f_i}$
for each $i \in \{1,\cdots,l\}$.
Hence each function $f_j$ is equal to
$\overline{f}_j \circ P$ for some
smooth $\overline{f}_j : V \to \R$ near $p$.
Since $\text{dim}(V) > l$
there is a smooth function
$\overline{f}_{l+1} : V \to \R$
equal to $0$ at $P(p)$ with
$d\overline{f}_{l+1}$
not in the span of $d\overline{f}_1,\cdots,d\overline{f}_l$ at $P(p)$.
Define $f_{l+1} := \overline{f}_{l+1} \circ P$.
Since $\iota_{X_{f_i}}df_{l+1} = 0$
for all $i =1,\cdots,l$ we get that
$f_{l+1}$ Poisson commutes with $f_1,\cdots,f_l$.
Hence by induction we have constructed
smooth functions
$f_1,\cdots,f_n$ defined near $p$ which all Poisson commute and so that $df_1,\cdots,df_n$
are linearly independent at each point
near $p$.
The result now follows from the
Darboux Theorem for fibrations
in \cite[Section 4.2]{dynamicalsystems4}
since the fibers of the map
$(f_1,\cdots,f_n)$ are Lagrangian near $p$.
\qed

\begin{lemma} \label{lemma:lagrangiansinstandardform}
Let $L \subset \widehat{M}$ be Lagrangian which is
transverse to $C_{I'}$ for all $I' \subset I$
and let $p \in C_I \cap L$ for some $I \subset \{1,\cdots,l\}$.
Then there exists a Lagrangian $L'$ so that:
\begin{enumerate}
\item $L'$ is Hamiltonian isotopic to $L$ where the Hamiltonian isotopy has support in an arbitrarily small neighborhood
$U$ of $p$,
\item $L'$ is transverse to $C_{I'}$ for all $I' \subset I$
and $\pi_I(L' \cap C_I) = \pi_I(L \cap C_I)$,
\item $L'$ is fiber radial near $p$.
\end{enumerate}
\end{lemma}
\proof
Choose a trivialization $\prod_{j \in I} A_j$ of
the fiber $\pi_I^{-1}(\pi_I(p))$
so that its real part $R$ passes through $p$.
We let $Q$ be a Lagrangian diffeomorphic to a ball so that:
\begin{itemize}
\item $Q \cap \pi_I^{-1}(\pi_I(p))$ is equal to
$R$ in a small neighborhood of $p$ and
\item $Q$ is transverse to $C_{I'}$ for all $I' \subset I$.
\end{itemize}

Let $B_\delta \subset \R^n$ be the open ball of radius $\delta$.
By Lemma \ref{lemma:darbouxchart}, we can construct
a Darboux chart centered at $p$
where some of the Darboux coordinates
are the commuting functions $(r_i - r_i(p))_{i \in I}$.
Using this chart we can then find
an open set $U$ containing $p$ so that:
\begin{itemize}
\item there is a symplectomorphism $\Phi : U \to V$ where
$V \subset T^* B_\delta$ is the set of covectors of norm $< \delta$ on $B_\delta$ for some $\delta < \frac{1}{4}$ and $\Phi(p)$ is the $0$ covector at $0 \in B_\delta$,
\item there is a submanifold $A \subset B_\delta$
so that $\Phi(C_I \cap U) = \pi_B^{-1}(A) \cap V$
where $\pi_B : T^* B_\delta \twoheadrightarrow B_\delta$ is the natural projection map,
\item for each $q \in V_I$ and each tuple $(t_i)_{i \in I}$ of constants,
$U \cap \cap_{i \in I}\{r_i = t_i\} \cap \pi_I^{-1}(q) \subset \Phi^{-1}(\pi_B^{-1}(b))$ for some $b \in B_\delta$
\item and $U$ is small enough so that there are smooth functions $f_1 : B_\delta \to \R$,
$f_2 : B_\delta \to \R$ so that
$\Phi(L \cap U)$ is equal to the image of the section $df_1$
and $\Phi(Q \cap U)$ is the image of the section $df_2$.
\end{itemize}

By subtracting appropriate constants we can assume that $f_1(\pi_B(\Phi(p))) = f_2(\pi_B(\Phi(p)))=0$.
Because $p \in L \cap Q$ and because $\Phi(p)$ is the $0$ covector at $0$, we get that $df_1 = df_2 = 0$ along $\pi_B(p)$.
Also because $R$ is transverse to $C_I$, we get that $A$ is transverse to the manifold $R' := \pi_B(\Phi(R \cap U))$ at $\pi_B(p)$ inside $B_\delta$.
Hence by Lemma \ref{lemma:lagrangiansinstandardform},
there is a smooth function $f : {\mathcal N} \to \R$ where ${\mathcal N}$ is an open set containing $\pi_B(p)$
so that $f|_{A \cap {\mathcal N}} = f_1|_{A \cap {\mathcal N}}$ and $f|_{R' \cap {\mathcal N}} = f_2|_{R' \cap {\mathcal N}}$
and that $df$ restricted to $R' \cap {\mathcal N}$
is equal to $df_2$ restricted to $R' \cap {\mathcal N}$.
Choose some $\delta' < \delta$ small enough so that
$B_{\delta'} \subset {\mathcal N}$,
and so that the norm of $d(f|_{B_{\delta'}})$
and $d(f_1|_{B_{\delta'}})$ is less than $\frac{1}{4}\delta$. 

Now let $\rho : B_{\delta'} \to \R$ be a smooth function with compact support
which is equal to $1$ near $\pi_B(p)$ and so that the norm of $d\rho$ is less than $2/{\delta'}$.
Define a smooth family of functions $(f^t : B_{\delta'} \to \R)_{t \in [0,1]}$ by
$f^t := (t (\rho f + (1-\rho)f_1) + (1-t) f_1)$.
Because $f(0)=f_1(0)$, $df(0) = df_1(0) =0$ and the norm of $d(f|_{B_{\delta'}})$ and $d(f_1|_{B_{\delta'}})$ is bounded above by $\frac{1}{4}\delta$, we get that
$|f - f_1| < \frac{1}{4}\delta \delta'$ inside $B_{\delta'}$.
This means that the norm of $df^t$ is less than $\delta$
and hence the graph of $df^t$ is contained in $V$.
Let $L^t$ be a smooth family of Lagrangians so that
$L^t = L$ outside $\Phi^{-1}(\pi_B^{-1}(B_{\delta'}))$
and $L^t$ is the graph of $df^t$ inside this set
(after taking its preimage under $\Phi$).
Then $L^t$ is induced by a Hamiltonian isotopy and the support of such a Hamiltonian
is contained in $U$ (which can be made arbitrarily small).
Hence $L' := L^1$ satisfies property $(1)$.
Also $L'$ is transverse to $C_{I'}$ for all $I' \subset I$
because the fibers of $\pi_B \circ \Phi : U \to B_\delta$ contain the level sets of the tuple
$(r_i)_{i \in I}$ inside $U$.
Because $f^t|_{A \cap B_{\delta'}} = f_1|_{A \cap B_{\delta'}}$
we get that $\pi_I(L' \cap C_I) = \pi_I(L \cap C_I)$.
Hence $(2)$ holds.
Finally because $df^1|_{R' \cap B_{\delta'}} = df_2|_{R' \cap B_{\delta'}}$,
we get that $L'$ is fiber radial near $p$ and so $(3)$ holds.
\qed

\proof of Lemma \ref{lemma:lagrangianpertubation}.

This Lemma is proven in three steps.
In step $1$, we use Lemma \ref{lemma:perturbinglagrangian}
and induction on the sets $I \subset \{1,\cdots,l\}$
to find a $C^\infty$ small Hamiltonians
$\phi_1,\phi_2$ whose support is disjoint from $M$
so that $\pi_I|_{\phi_1(L_1) \cap C_1}$
and $\pi_I|_{\phi_1(L_1) \cap C_2}$
intersect transversally
for all $I \subset \{1,\cdots,l\}$
and if $|I| < n$ then their
respective intersection points are isolated.
In other words, we show that a modified version of
(\ref{item:lemma441}) holds (note that we do not have a single Hamiltonian).
In step $2$ we find by Lemma
\ref{lemma:lagrangiansinstandardform}
$C^0$ small
Hamiltonian diffeomorphisms $\check{\phi}_1,\check{\phi}_2$
so that
$\check{\phi}_i(\phi_i(L_i))$ is fiber radial
near each point of
$\check{\phi}_i(\phi_i(L_i)) \cap \pi_I^{-1}(x) \cap C_I$ for all $i=1,2$
whenever $x$ satisfies
$\check{\phi}_1(\phi_1(L_1)) \cap \pi_I^{-1}(x) \cap C_I \neq 0$
and 
$\check{\phi}_1(\phi_2(L_2)) \cap \pi_I^{-1}(x) \cap C_I \neq 0$
while retaining the properties from Step 1.
In other words, a modified version of
both (\ref{item:lemma441}) and (\ref{item:lemma442}) holds.
Finally in Step 3, we construct our Hamiltonian
diffeomorphism $\phi$ using $\phi_1$ and $\phi_2$.

{\it Step 1}:
Let $C'_I := \cap_{j \in I} \{r_j = b^2+\epsilon_g\}$.
Then $C'_I$ is coisotropic and $C_I$ is an open subset of $C'_I$
whose closure is compact inside $C'_I$.
We choose a total ordering $\preccurlyeq$ on subsets $I \subset \{1,\cdots,k\}$ with $|I|<n$. We write $I \prec J$ if $I \preccurlyeq J$ and $I \neq J$.
We will induct on this ordering.
Suppose there are $C^\infty$ small
Hamiltonian diffeomorphisms
$\phi^\prec_i$, $i=1,2$ so that for all $I' \prec I$,
we have that $\pi_I|_{\phi^\prec(L_1) \cap C_I}$ and
$\pi_I|_{\phi^\prec(L_2) \cap C_I}$
 are transverse to each other and whose intersection points are isolated
on some neighborhood of the closure of $\pi_{I'}(C_{I'})$
and that these intersection points
are disjoint from $\partial \pi_{I'}(C_{I'})$.
We now wish to prove a similar statement for all $I' \preceq I$.
By using Lemma \ref{lemma:perturbinglagrangian} with $S = U_I,C = C'_I$, $B = \pi_I(C'_I)$ and with $K$ equal to the closure
of $\pi_I(C_I)$,
there are $C^\infty$ small Hamiltonian diffeomorphisms
$\phi_i^=$, $i=1,2$ with support near $K$,
so that $\pi_I|_{\phi^=(\phi^\prec(L_1)) \cap C_I}$,
$\pi_I|_{\phi^=(\phi^\prec(L_2)) \cap C_I}$ become transverse with isolated intersection points near $K$.
We can also ensure that these intersection points are disjoint from
$\partial \pi_{I}(C_{I})$.
We can make this perturbation small enough so that
the induction hypothesis still holds.
This completes the inductive step.
Hence,
we get $C^\infty$ small Hamiltonian
diffeomorphisms $\phi_1,\phi_2$ so that
$\pi_I|_{\phi_1(L_1) \cap C_1}$
and $\pi_I|_{\phi_2(L_1) \cap C_2}$
intersect transversally for all $I$
and so that they have isolated intersection points
whenever $|I| < n$.

{\it Step 2}:
By repeatedly applying
Lemma \ref{lemma:lagrangiansinstandardform}
to $\phi_1(L_1)$ and $\phi_2(L_2)$
near all points $p_1 \in \phi_1(L_1) \cap C_I$ and $p_2 \in \phi_2(L_2) \cap C_I$
where $\pi_I(p_1) = \pi_I(p_2)$ for some $I$,
we can find $C^0$ small Hamiltonian diffeomorphisms
$\check{\phi}_1,\check{\phi}_2$ so that
$\check{\phi}_1(\phi_1(L_1))$
and $\check{\phi}_2(\phi_2(L_2))$
become fiber radial near these points.
The support of these Hamiltonians can be made
disjoint from $C_{I'}$ for all $I'$ not contained in $I$ due to the fact that
$\pi_I(p_i)$ is disjoint from
$\partial \pi_I(C_I)$ for $i=1,2$.
We can also make sure the support of $\check{\phi}_1$
and $\check{\phi}_2$
is sufficiently small so that
for all $I' \subsetneq I$,
the immersions $\pi_{I'}|_{L_1 \cap C_{I'}}$
and
$\pi_{I'}|_{L_1 \cap C_{I'}}$
have no additional intersection points
and so that these immersions do not change near
the existing intersection points due to the fact that intersection points
$\pi_{I'}(L_1 \cap C_{I'}) \cap \pi_{I'}(L_2 \cap C_{I'})$ are disjoint from $\pi_{I'}(\partial C_{I'})$ which contains
$\pi_{I'}(C_I)$.
Also Lemma \ref{lemma:lagrangiansinstandardform}
tells us that $\pi_I|_{L_i \cap C_I}$ does not change.
Hence property $(1)$ still holds.

{\it Step 3}:
We can make sure that the support of
$\check{\phi}_i \circ \phi_i$
is contained in an arbitrarily small neighborhood
of $\cup_{I \subset \{1,\cdots,l\}} (C_I \cap L_i)$
for $i =1,2$.
In particular the supports of these Hamiltonian
diffeomorphisms can be made disjoint
because $L_1$ and $L_2$ are disjoint outside $M$
as they are transverse and cylindrical outside $M$.
We define $\phi$ to be the composition
$\check{\phi}_2 \circ \phi_2 \circ \check{\phi}_1 \circ \phi_1$.
Therefore $\phi(L_1)$,$\phi(L_2)$ satisfy
properties $(1)$ and $(2)$.
\qed

\bigskip

Let $\phi$ be the Hamiltonian symplectomorphism from Lemma \ref{lemma:lagrangianpertubation}.
We can choose $\phi$ so that
its support is disjoint from $M$
due to the fact that $\cup_i U_i$
is disjoint from $M$.
From now on we will replace our set of compatible annulus bundles $U_I,\pi_I,V_I$ on $\widehat{M}$
with its pullback by $\phi$.
In other words,
we replace $U_I$ with $\phi^{-1}(U_I)$ for all
$I$,
$\pi_I$ with $\pi_I \circ \phi$ for all
non-empty $I$
and we leave $V_I$ alone for all $I$.
This also means that we replace the coisotropic submanifolds $C_I$
with $\phi^{-1}(C_I)$.

In particular $L_1,L_2$ satisfy the following properties:
For each subset $I \subset \{1,\cdots,k\}$,
\begin{enumerate}
	\item \label{item:lemma4411}
	$L_1,L_2$ are transverse to $C_I$ for all $I$,
	the Lagrangian immersions $\iota_{1,I} := \pi_I|_{L_1 \cap C_I},\iota_{2,I} := \pi_I|_{L_2 \cap C_I}$ are transverse to each other
	and also the intersection points between $\iota_{1,I}$
	and $\iota_{2,I}$ are isolated if $|I|<n$,
	\item \label{item:lemma4422}
	and if
	$x \in \text{Image}(\iota_{1,I}) \cap \text{Image}(\iota_{2,I})$
	then $L_1$ (resp. $L_2$) is fiber radial near each point of
	$L_1 \cap \pi_I^{-1}(x) \cap C_I$
	(resp. $L_2 \cap \pi_I^{-1}(x) \cap C_I$).
\end{enumerate}

We will use Theorem \ref{theorem:growthratebound} in order to prove Theorem \ref{theorem:nicelystratifiedgrowth}.
Hence we need an appropriate Hamiltonian $H : \widehat{M} \to \R$
which is constructed as follows:
Define a smooth function $g : (b^2,B^2) \to (0,1)$ so that
\begin{gfun}
	\item $g(x)= 0$ for $x$ near $b^2$
	and $g(x) = (x-b^2)/ (B^2-b^2)$ for $x$ near $B^2$.
	\item \label{item:secondderivativeofg}
	We also assume that the derivative of $g$ is non-negative
	and that it is strictly positive when $g(x)$ is positive.
	We also assume that $g''$ is non-negative,
	and strictly positive only in the interval
	$((b+\epsilon_g)^2,c_g^2)$ for some $c_g$ satisfying $b+\epsilon_g < c_g<B$.
	\item We have that $g(x)$ is $0$ if and only if $x \in (b^2,(b + \epsilon_g)^2]$.
\end{gfun}

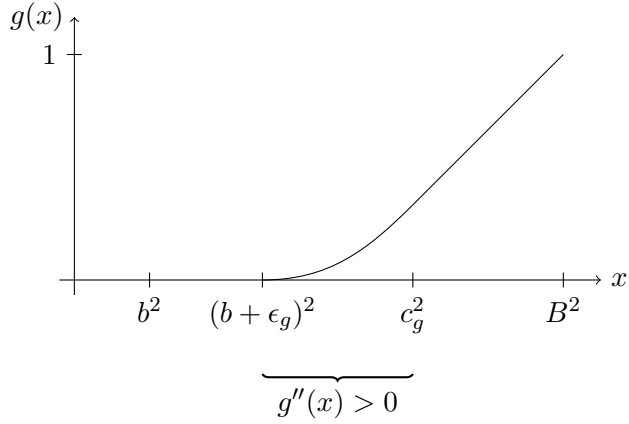
\begin{center}
	\begin{tikzpicture}[domain=0:4]
	\draw[->] (-1.2,0) -- (6,0) node[right] {$x$};
	\draw[->] (-1,-0.2) -- (-1,3.5) node[left] {$g(x)$};
	\draw (-0.9,3) -- (-1.1,3) node[left] {$1$};
	\draw (1.5,0) to[out=0,in=225] (3.5,1);
	\draw (3.5,1) -- (5.5,3);
	\draw (0,0.1) -- (0,-0.1) node[below,scale=1] {$b^2$};
	\draw (1.5,0.1) -- (1.5,-0.1) node[below,scale=1] {$(b+\epsilon_g)^2$};
	\draw (3.5,0.1) -- (3.5,-0.1) node[below,scale=1] {$c_g^2$};
	\draw (5.5,0.1) -- (5.5,-0.1) node[below,scale=1] {$B^2$};
	\draw [
	thick,
	decoration={
		brace,
		mirror,
		raise=0.5cm
	},
	decorate
	] (1.5,-0.75) -- (3.5,-0.75) 
	node [pos=0.5,anchor=north,yshift=-0.55cm] {$g''(x)>0$}; 
	\end{tikzpicture}
\end{center}

By abuse of notation we define $g(r_i^2)$
to be the function $U_\emptyset \to \R$
by defining it to be zero outside $U_i$
and $g(r_i^2)$ inside $U_i$.
By property \ref{item:closureproperty} combined with
the fact that $\{r_i > r\}$ is open for all $r \in (b, B)$,
we get that this function is smooth.
We define $K : U_\emptyset \to \R$ to be equal to $\sum_{i = 1}^l g(r_i^2)$.
We define $H : \widehat{M} \to \R$ to be a smooth function equal to $K$ near $\cap_i \{r_i \leq c_g\}$
and any strictly positive smooth function elsewhere which is constant at infinity.

\begin{lemma} \label{lemma:nondegenerateflow}
There is a small fixed neighborhood ${\mathcal N}$ of $H^{-1}(0)$
so that:
\begin{enumerate}
\item \label{item:transverselagrangians}
	$\phi^{\lambda H}_1(L_1) \cap {\mathcal N}$ and $L_2 \cap {\mathcal N}$
	are transverse to each other for all $\lambda \geq 0$.
\item \label{item:upperboundonfixedpoints}
	There is a polynomial $P$ of degree $C_D \leq \frac{1}{2}\text{dim}_{\R} \widehat{M}$
	so that for all $\delta >0$ sufficiently small, the number of elements in:
	\[\cup_{0 \leq \lambda' \leq \lambda} \phi^{\lambda' H}_1(L_1) \cap  L_2 \cap H^{-1}(\delta)\]
	and in:
	\[{\mathcal N} \cap \phi^{\lambda H}_1(L_1) \cap L_2\]
	is bounded above by $P(\lambda)$.
\end{enumerate}
\end{lemma}
\proof of Lemma \ref{lemma:nondegenerateflow}.
We will show that for all $I \subset \{1,\cdots,l\}$ and
$q \in H^{-1}(0) \cap C_I$
there is a neighborhood ${\mathcal N}_{q,I}$ of $q$ so that:
\begin{enumerate} [label=(\alph*)]
\item \label{item:transverselagrangiansnearq}
	$\phi^{\lambda H}(L_1) \cap {\mathcal N}_{q,I}$ and $L_2 \cap {\mathcal N}_{q,I}$
	are transverse to each other for all $\lambda \geq 0$,
\item \label{item:upperboundonfixedpointsnearq}
	there is a polynomial $P_{q,I}$ of degree $|I| \leq C_D$
	so that for all $\delta >0$ sufficiently small, the number of elements in:
	\[A_{q,I,\lambda} := \cup_{0 \leq \lambda' \leq \lambda} \phi^{\lambda' H}_1(L_1) \cap  L_2 \cap {\mathcal N}_{q,I} \cap H^{-1}(\delta)\]
	and in:
	\[B_{q,I,\lambda} := {\mathcal N}_{q,I} \cap \phi^{\lambda H}_1(L_1) \cap L_2\]
	is bounded above by $P_{q,I}(\lambda)$.
\end{enumerate}
This will prove the Lemma with ${\mathcal N}$ being a finite union of the neighborhoods ${\mathcal N}_{q,I}$
because $H^{-1}(0)$ is compact by \ref{item:closureproperty} and \ref{item:fiberproperty}.

Fix $q,I$ as above and let $F_q$ be the fiber of $\pi_I$ through $q$.
If $q \notin L_2$ then we choose ${\mathcal N}_{q,I}$ to be small enough so that
it is disjoint from $L_2$ and hence $\phi^{\lambda H}_1(L_1) \cap L_2 \cap {\mathcal N}_{q,I} = \emptyset$ for all $\lambda \geq 0$.
So from now on we will assume that $q \in L_2$.
Also if $L_1 \cap F_q = \emptyset$
then because the flow of $H$ preserves the fibers of $p_I$ and also the radial coordinates $(r_j)_{j \in I}$,
there is a neighborhood ${\mathcal N}_{q,I}$ of $q$ so that
$\phi^{\lambda H}_1(L_1) \cap L_2 \cap {\mathcal N}_{q,I} = \emptyset$ for all $\lambda \geq 0$.
Hence we will also assume that $L_1 \cap F_q \neq \emptyset$.

This means that $L_i$ is fiber radial
near each point of $L_i \cap C_I \cap F_q$ for $i=1,2$. Hence
there is a trivialization $\prod_{j \in I} A_j$
with associated polar coordinates $(r_j,\vartheta_j)_{j \in I}$
of $F_q$ so that:
\begin{itemize}
	\item there is a neighborhood ${\mathcal N}_{F_q,I}$
	of $C_I \cap F_q$ inside $\widehat{M}$ and constants
	$(t_j^k)_{j\in I,k = 1,\cdots,m}$ for some $m > 0$ so that $L_1 \cap F_q \cap {\mathcal N}_{F_q,I}$ is a disjoint union of $m$
	connected components:
	\[\sqcup_{k=1}^m \{\vartheta_j = t_j^k \, | \, \forall j \in I\}\cap F_q \cap {\mathcal N}_{F_q,I}.\]
	We also can assume that this neighborhood is invariant under the natural $U(1)^{|I|}$
	action on $U_I$.
	\item there is a neighborhood ${\mathcal N}'_{q,I}$ of $q$ so that
	$L_2 \cap F_q \cap {\mathcal N}'_{q,I}$ is equal to
	$L_2 = \{\vartheta_j = 0 \, | \, \forall j \in I\} \cap {\mathcal N}'_{q,I}$
\end{itemize}

Choose a small neighborhood ${\mathcal N}'_{F_q,I}$
of $F_q \cap C_I$ in $\widehat{M}$ invariant under the
$U(1)^{|I|}$ action on $U_I$ so that
\begin{equation} \label{eqn:intersectionpoint}
\pi_I(L_1 \cap {\mathcal N}'_{F_q,I}) \cap \pi_I(L_2 \cap {\mathcal N}'_{F_q,I}) = \{\pi_I(q)\}.
\end{equation}
Such a neighborhood exists because the intersection
points of $\iota_{1,I},\iota_{2,I}$ are discrete in $V_I$.
We define ${\mathcal N}_{q,I} := {\mathcal N}'_{q,I} \cap {\mathcal N}_{F_q,I} \cap {\mathcal N}'_{F_q,I} \cap \cap_{j \in I} \{r_j < c_g\}$
where $c_g$ is the constant in property \ref{item:secondderivativeofg} in the definition of $g$ above.

Because $\omega_M|_{F_q} = \sum_j \frac{1}{2} d(r_j^2) \wedge d\vartheta_j$ and $H|_{F_q} = \sum_{j \in I} g(r_j^2)$,
we have that $\phi^{\lambda H}_1|_{F_q} : F_q \to F_q$ satisfies:
\begin{equation} \label{eqn:hamiltonianflowformula}
\phi^{\lambda H}_1|_{F_q}\left((r_i,\vartheta_i)_{i \in I}\right) =
(r_i,\vartheta_i + 2 \lambda g'(r_j^2))_{i \in I}.
\end{equation}
Therefore $\phi^{\lambda H}_1(L_1) \cap F_q \cap {\mathcal N}_{F_q,I}$
is equal to:
\begin{equation} \label{eqn:flowformula}
  \sqcup_{k=1}^m \cap_{j \in I}\{\theta_j = t_j^k + 2 \lambda g'(r_j^2)\} \cap F_q \cap {\mathcal N}_{F_q,I}.
\end{equation}

Now that we have chosen ${\mathcal N}_{q,I}$,
we will show that property \ref{item:transverselagrangiansnearq} above holds.
Let $V \in T(L_1 \cap (\phi^{\lambda H}_1)^{-1}({\mathcal N}_{q,I}))$ be a non-zero vector at a point $p \in L_1 \cap (\phi^{\lambda H}_1)^{-1}({\mathcal N}_{q,I})$.
We wish to show that $D\phi^{\lambda H}_1(V) \notin TL_2 \cap {\mathcal N}_{q,I}$ for all $\lambda \geq 0$.
If $V \notin T\widehat{M}|_{F_q}$ then $\phi^{\lambda H}_1(p) \notin L_2 \cap {\mathcal N}_{q,I}$
for all $\lambda \geq 0$ and so $D\phi^{\lambda H}_1(V) \notin TL_2$ for all $\lambda \geq 0$.
If $V \in T\widehat{M}|_{F_q}$ then $V = \sum_j a_j \frac{\partial}{\partial r_j} + Z$
where $Z$ is symplectically orthogonal to $F_q$ for some constants $(a_j)_{j \in I}$ due to the fact that $L_1$
is fiber radial near each point of $L_1 \cap F_q \cap {\mathcal N}_{F_q,I}$.
Hence by Equation (\ref{eqn:hamiltonianflowformula}),
$D\phi^{\lambda H}_1 (V) = \sum_{j \in I} 4 a_i r_i g''(r_i^2) \frac{\partial}{\partial \vartheta_i} + V$
at the point $\phi^{\lambda H}_1(p)$.
This cannot be tangent to $TL_2$ for any $\lambda$
because $\phi_1^{\lambda H}(p) \in {\mathcal N}_{q,I}$,
$g'(r_i),g''(r_i)>0$ (as $b + \epsilon_g < r_i < c_g$)
and because any vector tangent to $TL_2$ at any
 point in $F_q \cap {\mathcal N}'_{q,I}$
is equal
$\sum_{j \in I} b_j \frac{\partial}{\partial r_j} + W$
for some constants $(b_j)_{j \in I}$
where $W$ is symplectically orthogonal to $F_q$.
Hence \ref{item:transverselagrangiansnearq} holds.

We will now show that
property \ref{item:upperboundonfixedpointsnearq} holds.
This is just done by examining Equation (\ref{eqn:flowformula}) as every intersection point
$\phi^{\lambda H}_1(L_1) \cap L_2 \cap {\mathcal N}_{q,I}$ is contained in $F_q$ by
Equation (\ref{eqn:intersectionpoint}).
This just means we only need to find out when the angle coordinates $(\vartheta_j)_{j \in I}$
vanish inside $F_q$ when we flow from a point in $(\phi^{\lambda H}_1)^{-1}({\mathcal N}_{q,I} \cap F_q)$.

For each $1 \leq k \leq m$, let $\kappa^k : (0,\infty) \times (b+\epsilon_g, c_g)^{|I|} \to \R^{|I|}$
send $(\lambda, (x_j)_{j \in I})$ to
$(t_j^k + 2 \lambda g'(r_j^2))_{j \in I}$.
Let $\nu_\delta \subset (b+\epsilon_g, c_g)^{|I|}$
be equal to the set 
\[\{ (x_i)_{j \in I} \in (b+\epsilon_g,c_g)^{|I|} \, | \, \sum_{j \in I} g(x_j) = \delta \}.\]
For each $1 \leq k \leq m$, let
$\kappa^k_{1,\lambda} := \kappa|_{(0,\lambda) \times \nu_\delta}$
and
$\kappa^k_{2,\lambda} := \kappa|_{\{\lambda\} \times (b+\epsilon_g,c_g)^{|I|}}$.
Equation (\ref{eqn:flowformula})
then tells us that the number of points in $A_{q,I,\lambda}$
 is bounded above by the number of points in
$\cup_{k=1}^m (\kappa^k_{1,\lambda})^{-1}( (2\pi \Z)^{|I|})$
and the number of points in $B_{q,I,\lambda}$ is bounded above by the number of points in
$\cup_{k=1}^m (\kappa^k_{2,\lambda})^{-1}( (2\pi \Z)^{|I|})$.
Because $g',g''>0$ inside $(b+\epsilon_g,c_g)$,
we get that $\kappa^k_{1,\lambda}$
and $\kappa^k_{2,\lambda}$ are injective maps.
Also their images are contained a ball of radius proportional to $\lambda$.
This implies that the number of points in
$(\kappa^k_{1,\lambda})^{-1}( (2\pi \Z)^{|I|})$
and in
$(\kappa^k_{2,\lambda})^{-1}( (2\pi \Z)^{|I|})$
is bounded above by $C \lambda^{|I|}$ for some constant $C>0$ independent of $k$ and $\lambda$.
Hence the number of points in
$A_{q,I,\lambda}$ and $B_{q,I,\lambda}$
is bounded above by $m C \lambda^{|I|}$
which is a polynomial of degree $|I|$ in $\lambda$.
\qed


\bigskip

\proof of Theorem \ref{theorem:nicelystratifiedgrowth}.
We have two constants $0 < \delta_H \ll C_H \ll 1$ so that:
\begin{PB}

\item $H^{-1}((-\infty,C_H \pm \delta_H])$ is compact and contains $M$ by \ref{item:closureproperty} and \ref{item:fiberproperty}
and \ref{item:secondderivativeofg}.
\item For all $C \in [C_H - \delta_H,C_H+\delta_H]$,
$C$ is a regular value of $H$ and
$L_1$ and $L_2$ intersect $H^{-1}(C)$ transversally
by \ref{item:secondderivativeofg}.
\item 
By Lemma \ref{lemma:nondegenerateflow},
 we have for all $\lambda \geq 0$ that
$\phi^H_\lambda(L_1)$ and $L_2$ intersect transversally inside a small neighborhood of $H^{-1}([C_H-\delta_H,C_H+\delta_H])$.
Also the number of such intersection points is bounded above by $P(\lambda)$
where $P$ is a polynomial of degree $C_D$.
\item 
Lemma \ref{lemma:nondegenerateflow}, we have
for all $C \in [C_H - \delta_H,C_H+\delta_H]$ that
the number of flowlines of  $X_H$ inside $H^{-1}(C)$ of length $\leq \lambda$ starting on $L_1$ and ending on $L_2$
is bounded above by $P(\lambda)$.
\end{PB}

Hence $H$ is $(L_1,L_2,P)$-bounded where $P$ is a polynomial of degree $C_D$.
So by Theorem \ref{theorem:growthratebound} we get that $\Gamma(L_1,L_2;\K) \leq C_D$.
\qed

\section{Proof of the main theorem} \label{section:proofofmaintheorem}

We start with a Liouville domain $(M,\theta_M)$.
Here is the statement of Theorem \ref{theorem:main}:
{\it Suppose that $\partial M$ is contactomorphic to
	the link of an isolated complex singularity
	or $\widehat{M}$ is symplectomorphic to a
	smooth affine variety.
	Then $\Gamma(L_1,L_2;\K) \leq n$  for any
	transversally intersecting admissible Lagrangians $L_1,L_2$ in $\widehat{M}$ and any field $\K$ where $n$ is the complex dimension of
	our variety.
	}

This theorem follows immediately from
Theorem \ref{theorem:nicelystratifiedgrowth}
combined with Propositions
\ref{prop:isolatedsingularitymaintheorem}
and \ref{prop:smoothaffinemaintheorem} below.

\subsection{Proof of the Main Theorem for Isolated Complex Singularities}

Let $B \subset \C^N$
be an affine variety which has an isolated
complex singularity at $0$.
The {\bf link} $L_B$ of $B$
is defined to be $B \cap S_\epsilon$
where $S_\epsilon$ is a sphere of radius $\epsilon>0$ in $\C^N$.
For $\epsilon$ small enough $L_B$ has a natural contact structure
$\xi_B$ given by $TL_B \cap J_{\text{std}} TL_B$
where $J_{\text{std}}$ is the standard complex structure on $\C^N$
(see \cite{Varchenko:isolatedsingularities}).

Suppose we have a resolution $\pi : \widetilde{B} \to B$
of our isolated singularity $B$ by blowing it up along smooth loci so that $\pi^{-1}(0)$ is a union of transversally intersecting
complex hypersurfaces $E_1,\cdots,E_l$.
Such a resolution exists by \cite{hironaka:resolution}.
Let $C_\pi$ be size of the largest set $I \subset \{1,\cdots,l\}$
so that $\cap_{j \in I} E_j \neq \emptyset$.

\begin{prop} \label{prop:isolatedsingularitymaintheorem}
If $\partial{M}$ is contactomorphic to $L_B$
then $\widehat{M}$ has admits compatible annulus bundles at infinity of codimension $C_\pi$.
\end{prop}

Such a contactomorphism may be coorientation reversing.
We begin with the following definition:
\begin{defn} \label{defn:niceneighborhood}
Let $(X,\omega)$ be a symplectic manifold and let
$(S_1,\cdots,S_l)$ be codimension $2$ transversally intersecting symplectic submanifolds
so that $S_I := \cap_{j \in I} S_j$ are also symplectic
for all $I \subset \{1,\cdots,l\}$.
A {\bf compatible disk bundle neighborhood} of $\cup_i S_i$ consists of
neighborhoods $U_I$ of $S_I$ and fibrations $\pi_I : U_I \twoheadrightarrow S_I$
satisfying the following requirements:
\begin{enumerate}
\item $U_I \cap U_J = U_{I \cup J} \quad \forall I,J \subset \{1,\cdots, l\}$.
\item \label{item:propertyehresmannstructure}
The fibration $\pi_I$ has a $U(1)^{|I|}$ structure
group whose fibers are equal to a product
$\prod_{j \in I} \D_j$
where $\D_j \subset \C$ is the disk of radius $\epsilon$
with the standard symplectic form
and where the $U(1)$ factor corresponding to $i \in I$ rotates the disk $\D_i \subset \C$.
The natural Ehresmann connection consisting of vectors symplectically orthogonal to the fibers
of $\pi_I$ is compatible with this $U(1)^{|I|}$ structure group.
\item \label{item:propertyfiberform}
 For each $I \subset J$, we have that the fibers of $\pi_I|_{U_J}$ are contained in the fibers of $\pi_J$,
and if we have a trivialization $\prod_{j \in J} \D_j$ of a fiber $F$ of $\pi_J$ then
the fibers of $\pi_I|_F$ are of the form:
$\prod_{j \in I} \D_j \times \prod_{j \in I \setminus J} \{z_j\}$.
\end{enumerate}
\end{defn}
Define $U_i := U_{\{i\}}$
and $\pi_i := \pi_{\{i\}}$.
We have a natural function $r_i : U_i \to \R$
whose restriction to each trivialized fiber $\D_i$ of
our disk fibration $\pi_i$ is the natural radial
coordinate on $\D_i$ which we call the {\bf $i$th radial coordinate}.

\begin{lemma} \label{lemma:niceneighborhoodimpliesnicelystratified}
Let $(W,\omega)$ be a symplectic manifold which is connected
and let
$(S_1,\cdots,S_l)$ be compact codimension $2$ transversally intersecting symplectic submanifolds
so that $S_I := \cap_{j \in I} S_j$ are also symplectic
for all $I \subset \{1,\cdots,l\}$
and suppose that $\cup_i S_i$ admits a compatible disk bundle neighborhood.
Let $C_S$ be the cardinality of the largest set $I$ satisfying $S_I \neq \emptyset$.
Suppose that $\omega|_{W \setminus \cup_i S_i} = d\theta$ for some $1$-form $\theta \in \Omega^1(W \setminus \cup_i S_i)$
and let $f : W \setminus \cup_i S_i \to \R$ be a smooth function with proper level sets
such that $f$ tends to $-\infty$ as we approach $\cup_i S_i$
and $df(X_\theta) \neq 0$ at every point in $W \setminus \cup_i S_i$.
Then any Liouville domain $(M,\theta_M)$ with connected boundary contactomorphic to
$(f^{-1}(c),\theta|_{f^{-1}(c)})$ for some $c$
admits compatible annulus bundles at infinity of codimension $C_S$ as in Definition \ref{definition:nicelystratifiedatinfinity}.
\end{lemma}
\proof of Lemma \ref{lemma:niceneighborhoodimpliesnicelystratified}:
Let $U_I,\D_j,\pi_I,\epsilon$ be our compatible disk bundle neighborhood
as in Definition  \ref{defn:niceneighborhood}.
Let $r_i$ be the $i$th radial coordinate.
We suppose that each $S_i$ is connected.
We let $(r_i,\vartheta_i)$ be the standard polar coordinates
for $\D_i$.
For a subset $E \subset [0,\epsilon]$ and each
nonempty $I \subset \{1,\cdots,l\}$, define:
$U^{E}_I \subset U_I$
to be $U_I \cap \cap_{i \in I} \{r_i^{-1}(E)\}$.
Define
\[\dot{U}_I := U^{(\frac{1}{4}\epsilon,\frac{3}{4}\epsilon)}_I \setminus \cup_{j \notin I} U^{[0,\frac{1}{4}\epsilon]}_{j}.\]
Define $\dot{U}_i := \dot{U}_{\{i\}}$
and $\dot{U} := \cup_{i=1}^l \dot{U}_i$.


Because $\dot{U}$ is relatively compact,
there are constants $c_1 < c < c_2$
so that
$\dot{U} \subset f^{-1}((c_1,c_2))$.
There is a constant $C > 0$
with the property that
the length of every flowline
of $X_{\theta}$ or $-X_{\theta}$
from a point in $f^{-1}(c)$
to a point in $f^{-1}(\{c_1,c_2\})$
is at most $C$.

Let
$\Phi : f^{-1}(c) \to \partial M$
be our contactomorphism.
There is a smooth function
$g : f^{-1}(c) \to \R \setminus \{0\}$
with the property that
$\Phi^*(\alpha_M) = g \theta|_{f^{-1}(c)}$
where $\alpha_M := \theta_M|_{\partial M}$.
Here $g>0$ if and only if $\Phi$ is coorientation preserving.
If $g > 0$, we let $m > 0$ be smaller than the infimum of $1/g$.
If $g < 0$, we let $m < 0$ be larger than the supremum of $1/g$.

We then have an embedding
$\iota_U : f^{-1}([c_1,c_2]) \hookrightarrow \widehat{M}$
into the cylindrical end defined as follows:
Let $p_c : f^{-1}([c_1,c_2]) \to f^{-1}(c)$
send $x$ to the unique intersection point
of the flowline of $X_\theta$
passing through $x$ with $f^{-1}(c)$.
Let $\phi^x : (m_x,M_x) \to W \setminus \cup_i S_i$
be the unique largest flowline of $X_\theta$
so that $\phi^x(0) = x$ where $m_x < 0$ and $M_x > 0$.
Define $l : f^{-1}([c_1,c_2]) \to \R$
by $l(x) := t$ where $t \in (m_x,M_x)$
is the unique point satisfying
$\phi^x(t) \in f^{-1}(c)$.
We define:
$\iota_U(x) := \left(\frac{e^{-l(x)}e^C}{g(p_c(x)) m}, \Phi(p_c(x))\right)
\in [1,\infty) \times \partial M \subset \widehat{M}$.
This satisfies
$\iota_U^*(d\theta_M) = (e^C/m) \omega$
due to the fact that the $\omega_M$
dual $X_{\theta_M}$ of $\theta_M$
is $r_M \frac{\partial}{\partial r_M}$ inside
the cylindrical end of $M$.
Also
$(\iota_U)_*(X_\theta|_{f^{-1}([c_1,c_2])}) = X_{\theta_M}|_{\text{Image}(\iota_U)}.$

If $df(X_\theta)>0$,
define
$B := \iota_U(f^{-1}(c_1))$.
Otherwise define
$B := \iota_U(f^{-1}(c_2))$.
In other words, we are choosing $B$
so that $X_{\theta_M}$ points inwards along $B \subset \iota_U(f^{-1}([c_1,c_2]))$.
If $df(X_\theta)>0$,
define
$U_B := \iota_U(\dot{U} \cup (f^{-1}([c_1,c_2]) \cap \cup_i U_i^{[0,\frac{1}{4}\epsilon])})$,
otherwise define
$U_B := \iota_U(\dot{U} \cup (f^{-1}([c_1,c_2]) \setminus \cup_i U_i^{(0,\frac{3}{4}\epsilon)}))$.
Morally, $U_B$ is the region which
`fills' in $\iota_U(\dot{U}) \subset \iota_U(f^{-1}([c_1,c_2]))$.
It will enable us to define $U_\emptyset$.

Because $B$ is isotopic in $\widehat{M}$
to $\partial M$ through smooth embedded hypersurfaces,
we get that $B$ is the boundary of a unique compact
codimension $0$ submanifold $M_B$
containing $M$.
We define $\check{U}_\emptyset := U_B \cup M_B$.
Define $\check{V}_\emptyset := U_\emptyset$
and $\check{\pi}_\emptyset : \check{U}_\emptyset \to \check{V}_\emptyset$
to be the identity map.
For all non-empty $I \subset \{1,\cdots,l\}$,
define:
$\check{U}_I := \iota_U(\dot{U}_I)$,
$\check{V}_I := \pi_I(\dot{U}_I)$
and
$\check{\pi}_I : \check{U}_I \to \check{V}_I$
by
$\check{\pi}_I(x) := \pi_I(\iota_U^{-1}(x))$.
Define $\check{U}_i := \check{U}_{\{i\}}$.

Then: $\check{U}_{I \cup J} = \check{U}_I \cap \check{U}_J$
for all $I,J \subset \{1,\cdots,l\}$.
Because the image of $\iota_U$
is disjoint from $M$, we get that
$\check{U}_\emptyset \setminus \cup_i \check{U}_{i}$
contains $M$.
If $|I| > C_S$ then $\check{U}_I = \emptyset$
because $U_I = \emptyset$.
Hence \ref{item:containsa} is satisfied.

Property (\ref{item:propertyehresmannstructure})
 of Definition
\ref{defn:niceneighborhood} implies that
\ref{item:structuregroupcondition} is satisfied.
Because $M_B$ and $U_B$ are relatively compact,
we get that $U_\emptyset$ is relatively compact.
Hence $\check{U}_\emptyset \setminus \cup_i \check{U}_{i}$
is relatively compact.
Also for non-empty $I \subset\{1,\cdots,l\}$,
because $\pi_I(\dot{U}_I \setminus \cup_{j \notin I} \dot{U}_{j}) = \{S_I \setminus \left( \cup_{j \notin I}\{r_j < \frac{3}{4}\epsilon\}\right)$ is relatively compact
inside $\pi_I(\dot{U}_I) = \{S_I \setminus \left( \cup_{j \notin I}\{r_j \leq \frac{\epsilon}{4}\}\right)$,
we get that $\pi_I(\check{U}_I \setminus \cup_j \check{U}_{j})$
is relatively compact inside
$\check{V}_I$.
These facts combined with
part (\ref{item:propertyfiberform})
of Definition \ref{defn:niceneighborhood}
show that
\ref{item:fiberproperty} is satisfied.

Define $r_i^* : \dot{U}_i \to \R$
to be $r_i$ if $df(X_\theta) >0$
and $\epsilon -r_i$ otherwise.
The functions
$\check{r}_i : \check{U}_{i} \to \R$
defined by
$\check{r}_i(x) := r_i^*(\iota_U^{-1}(x))$
are equal to the natural radial coordinates
on the annulus fibers of $\pi_{i}$
for each $i \in \{1,\cdots,l\}$.
Because
$\{r^*_i \geq r\} \subset f^{-1}((c_1,c_2))$
is closed inside $\iota_U^{-1}(U_B)$ and because
$U_B \cap M = B$,
we get that
$\{\check{r}_i \geq r\}$ is closed inside
$U_\emptyset$.
Hence \ref{item:closureproperty}
is satisfied.

\qed

\proof of Proposition \ref{prop:isolatedsingularitymaintheorem}.
By applying \cite[Theorem 5.25]{McLean:isolated}
to our chosen resolution $\pi$
combined with \cite[Theorem 2.12]{McLeanTehraniZinger:smoothing}
(or equivalently
\cite[Theorem 5.3, Theorem 5.20]{McLean:affinegrowth})
and \cite[Corollary 5.11]{McLean:isolated},
there is a compact symplectic manifold $(X,\omega)$
(which is a codimension $0$ submanifold of our resolution),
codimension $2$ closed submanifolds $S_0, \cdots, S_l \subset X$ (corresponding to the exceptional divisors of this resolution),
a $1$-form $\theta$ on $X \setminus \cup_i S_i$
and a smooth function $f \in C^\infty(X \setminus \cup_i S_i)$ so that:

\begin{enumerate}
\item $f$ tends to $-\infty$ as we approach $\cup_i S_i$,
and $f^{-1}(c)$ is compact for all sufficiently negative $c$.
\item $\omega|_{X \setminus \cup_i S_i} = d\theta$ and
$df(X_\theta)>0$ along sufficiently negative level sets of $f$
and $(f^{-1}(c),\theta|_{f^{-1}(c)})$
is a contact manifold
contactomorphic to $(L_B,\xi_B)$ for all sufficiently negative $c$.
\item If $|I| > C_\pi$ then $S_I = \emptyset$.
\item $S_I$ is a symplectic submanifold for all
$I \subset \{1,\cdots,l\}$ and
$\cup_i S_i$ admits a compatible disk bundle neighborhood.
\end{enumerate}
Then by Lemma \ref{lemma:niceneighborhoodimpliesnicelystratified} with
$W = \{f < -C\} \cup \cup_i S_i$
for some sufficiently large $C$ we get that
$\widehat{M}$ admits compatible annulus bundles at infinity of codimension $C_\pi$.
\qed

\subsection{Proof of the Main Theorem for Smooth Affine Varieties}

Let $A \subset \C^N$ be a smooth affine variety with
a symplectic structure $\omega_A$ given
by restricting the standard symplectic structure
on $\C^N$ to $A$.
We can view it as an open subset of some projective variety.
By \cite{hironaka:resolution} we can blow up this projective variety away from $A$ so that $A$
becomes an open subset of a smooth projective variety $X$
where $D := X \setminus A$ is a smooth normal crossing divisor.
Let $C_D \in \N$ be the codimension of the strata of $D$ of lowest dimension.

\begin{prop} \label{prop:smoothaffinemaintheorem}
If $\widehat{M}$ is symplectomorphic to $A$
then $\widehat{M}$ admits compatible annulus bundles at infinity of codimension $C_D$.
\end{prop}

We need some definitions and lemmas before we prove this theorem.
\begin{defn} (\cite[1.7.1]{EliahbergGromov:convexsymplecticmanifolds}, \cite[Section 2]{SS:rama}).
A {\bf convex symplectic manifold} is a manifold $N$ with a $1$-form $\theta_N$
such that
\begin{enumerate}
\item $\omega_N := d\theta_N$ is a symplectic form.
\item There is an exhausting function $f_N : N \rightarrow \R$
and a sequence $c_1 < c_2 < \cdots$ tending to infinity such that
the $\omega_N$-dual $X_{\theta_N}$ of $\theta_N$ satisfies
$df_N(X_{\theta_N}) > 0$ along $f_N^{-1}(c_i)$ for each $i$.
(Recall that an {\bf exhausting function}
is a smooth function which is proper and bounded from below.)
\end{enumerate}
We say that $(N,\theta_N)$ is of {\bf finite type} if there is an exhausting function $f_N : N \to \R$
and a constant $C \in \R$
such that $df_N(X_{\theta_N})>0$ along $f_N^{-1}(c)$
for all $c \geq C$.
\end{defn}

One important example of a finite type convex symplectic manifold is the completion of a Liouville domain
$\widehat{M}$ where $f_{\widehat{M}}$ is an exhausting function equal to $r_M$ at outside a compact set.

\begin{defn}
Let $(N,\theta^t_N)$ be a smooth family of convex symplectic manifolds
parameterized by $t \in [0,1]$.
This is said to be a {\bf convex deformation equivalence} if
for every $t \in [0,1]$ there is a
constant $\delta_t>0$,
an exhausting function $f_N^t : N \to \R$ and a sequence
of constants $c_1^t < c_2^t < \cdots$ tending to infinity
such that $df_N^t(X_{\theta^s_N}) > 0$
along $(f_N^t)^{-1}(c_i^t)$
for each $s \in [t- \delta_t,t+\delta_t]$ and each $i \in \N$.
We do not require that $f_N^t$,$c_i^t$,$\delta_t$
smoothly varies with $t$. In fact it can vary in
a discontinuous way with $t$.
\end{defn}

\bigskip

An affine variety $A \subset \C^N$ has a $1$-form $\theta_A := \sum_j \frac{1}{2}r_j^2 d\vartheta_j|_A$ where $(\vartheta_j,r_j)$ are polar coordinates for the $j$th $\C$ factor.
By \cite[Lemma 2.1]{McLean:affinegrowth} we have that $(A,\omega_A)$ is symplectomorphic to $(\widehat{\overline{A}},\theta_A)$
where $(\overline{A},
\theta_A)$ is the Liouville domain given by intersecting $A$ with a large closed ball in $\C^N$.

\proof of Proposition \ref{prop:smoothaffinemaintheorem}.

The smooth projective variety $X$ admits a Fubini-Study symplectic form.
By \cite[Theorem 5.20]{McLean:affinegrowth} there are compact codimension $2$ symplectic submanifolds $S_1,\cdots,S_l$ of $X$ so that
\begin{itemize}
	\item for every
	$I \subset \{1,\cdots,l\}$, we have that $\cap_{i \in I} S_i$
	is symplectic,
	\item  $\cup_i S_i$
	admits a compatible disk bundle neighborhood
	of codimension $C_D$
	\item  and so that $X \setminus \cup_i S_i$
	has the structure of a finite type convex symplectic manifold $(M_1,\theta_1)$
	convex deformation equivalent to $A = \widehat{\overline{A}}$.
\end{itemize} 
Because $(M_1,\theta_1)$ is a finite type convex symplectic manifold,
there exists an exhausting function $f_1 : M_1 \to \R$ and a constant $C \in \R$
so that $df_1(X_{\theta_1}) > 0$ along $f_1^{-1}(c)$ for all $c \geq C$.
Define $\overline{M}_1 := f_1^{-1}((-\infty,C])$.
Then $(\overline{M}_1,\theta_1)$ is a Liouville domain.
Define $f := -f_1$.
By Lemma \ref{lemma:niceneighborhoodimpliesnicelystratified},
with the symplectic manifold $W = f^{-1}(-\infty,-C) \cup \cup_i S_i \subset X$ and the function $f|_{W \setminus \cup_i S_i}$ and
the $1$-form $\theta_1|_{W - \cup_i S_i}$,
we get that $\widehat{\overline{M}}_1$ admits compatible annulus bundles at infinity
$U_I,V_I,\pi_I$
of codimension $C_D$
as in Definition \ref{definition:nicelystratifiedatinfinity}.
Because $M_1$ is convex deformation
equivalent to $\widehat{\overline{M}}_1$
by \cite[Corollary 8.3]{McLean:affinegrowth},
and because it is also convex deformation
equivalent to
$A = \widehat{\overline{A}}$, we get
that $\widehat{\overline{M}}_1$ is convex deformation
equivalent to $\widehat{\overline{A}}$.
So by \cite[Corollary 8.6]{McLean:affinegrowth},
$\widehat{M}_1$ and $\widehat{\overline{A}}$
are symplectomorphic.
Hence we have
a symplectomorphism
$\Phi : \widehat{M} \to \widehat{\overline{M}}_1$.

Let $\phi_t : \widehat{\overline{M}}_1 \to \widehat{\overline{M}}_1$
be the flow of the
Liouville vector field $X_{\theta_1}$.
Because $U_I \subset M_1$
is disjoint from $\overline{M}_1$ by \ref{item:containsa}
and is relatively compact for all $I$
by \ref{item:fiberproperty},
there exists
$T>0$ so that $\Phi(M) \subset \phi_T(U_I)$
for all $I \neq \emptyset$.
Hence $U'_I := \Phi^{-1}(\phi_T(U_I)), V'_I := V_I,
\pi'_I := \pi_I \circ (\phi_T)^{-1} \circ \Phi|_{U'_I}$
are compatible annulus bundles at infinity
of codimension $C_D$ on $\widehat{M}$.
\qed

\bibliography{references}

\end{document}